\documentclass[12 pt, draftcls, onecolumn]{IEEEtran}
\usepackage{amssymb,epsfig,color,cite,amsmath,amsfonts,mathrsfs,algorithm,algorithmic}
\usepackage{epstopdf}

\newtheorem{lemma}{Lemma}[section]
\newtheorem{theorem}{Theorem}[section]

\newtheorem{remark}{Remark}[section]
\newtheorem{definition}{Definition}[section]

\newtheorem{proposition}{Proposition}[section]

\allowdisplaybreaks

\begin{document}
%
\title{Distributed Continuous-Time Algorithm for Constrained Convex Optimizations via Nonsmooth Analysis Approach}
\title{
Distributed Computation of Linear Matrix Equations: An Optimization Perspective
}

\author{Xianlin~Zeng,
        Shu~Liang,
        Yiguang~Hong, and Jie Chen
\thanks{X. Zeng (xianlin.zeng@bit.edu.cn) is with the Key Laboratory of Intelligent Control and Decision of Complex Systems, School of Automation, Beijing Institute of
Technology, 100081, Beijing, China.}
\thanks{S. Liang (sliang@amss.ac.cn) is with the Key Laboratory of Knowledge Automation for Industrial Processes of Ministry of Education, School of Automation and
Electrical Engineering, University of Science and Technology Beijing, 100083, Beijing,  China.}
\thanks{Y. Hong (yghong@iss.ac.cn) is with the Key Laboratory of Systems and Control,  Academy of Mathematics and Systems Science, Chinese Academy of Sciences,  100190, Beijing, China. }
\thanks{J. Chen (chenjie@bit.edu.cn) is with Beijing Advanced Innovation Center for Intelligent Robots and Systems (Beijing Institute of Technology), Key Laboratory of Biomimetic Robots and Systems (Beijing Institute of Technology), Ministry of Education, 100081, Beijing,  China.}
}


\maketitle

\begin{abstract}

This paper investigates the distributed computation of the well-known linear matrix equation in the form of $AXB = F$, with the matrices $A$, $B$, $X$, and $F$ of appropriate dimensions, over multi-agent networks from an optimization perspective. In this paper, we consider the standard distributed matrix-information structures, where each agent of the considered multi-agent network has access to one of the sub-block matrices of $A$, $B$, and $F$. To be specific, we first propose different decomposition methods to reformulate the matrix equations in standard structures as distributed constrained optimization problems by introducing substitutional variables; we show that the solutions of the reformulated distributed optimization problems are equivalent to least squares solutions to original matrix equations; and we design distributed continuous-time algorithms for the constrained optimization problems, even by using augmented matrices and a derivative feedback technique. Moreover, we prove the exponential convergence of the algorithms to a least squares solution to the matrix equation for any initial condition.
\end{abstract}


%

\section{Introduction}

Recently, the increasing scale and big data of engineering systems and science problems have posed new challenges for the design based on computation, communication, and control. Traditional centralized algorithms  for the computation of small or modest sized problems are often entirely infeasible for large-scale problems.
As a result, distributed algorithms over multi-agent networks have attracted a significant amount of research attention due to their broad range of applications in nature science, social science and engineering.
Particularly, distributed optimization, which seeks a global optimal solution with the objective function as a sum of the local objective functions of agents, has become increasingly popular \cite{NOP:2010,YHL:SCL:2015,KCM:A:2015}.
In fact, distributed optimization with different types of constraints, including local constraints and coupled constraints, has been considered and investigated using either discrete-time  or continuous-time  solvers (see \cite{NOP:2010,YHL:SCL:2015,KCM:A:2015,QLX:A:2016,YHX:Cybernetics:2016,SJ:2012,ZYH:2017}).
Recently, distributed continuous-time algorithms have received much attention in \cite{WE:2011,GC:2014,YHL:SCL:2015,KCM:A:2015,QLX:A:2016,ZYH:2017,LW:TAC:2015,SA:TAC:2015,MLM:TAC:2015}, mainly because the continuous-time physical system may involve with solving optimal solutions and the continuous-time approach may provide an effective tool for analysis and design, though  distributed designs for many important problems  are still challenging.

In fact, distributed computation of the linear algebraic equation of the form $Ax=b$, where $A$ is a matrix and $x$ and $b$ are vectors of appropriate dimensions, over a multi-agent network  has attracted much research attention, because it is fundamental for many computational tasks and practical engineering problems.  Mainly based on the distributed optimization idea, distributed algorithms appeared for solving the  linear algebraic equation $Ax=b$.   The significant results in \cite{SA:TAC:2015,WFM:ACC:2016,SA:ACC:2016,MLM:TAC:2015, AMMH:Arxiv:2015,LLAS:Arxiv:2017,LMM:TAC:2017,LMNM:Automatica:2017} provided various distributed algorithms with the standard case that each agent knows a few rows of $A$ and $b$, while \cite{CZH:CCC:2017} proposed a distributed computation approach for another standard case, where each agent has the knowledge of a few columns of matrix $A$. 
In fact, the analysis given in \cite{MLM:TAC:2015,AMMH:Arxiv:2015,CZH:CCC:2017,LMM:TAC:2017,LMNM:Automatica:2017} depend on the existence of  exact solutions to the linear equations.
Specifically, \cite{MLM:TAC:2015} proposed a discrete-time distributed algorithm for a solvable linear equation and presented the necessary and sufficient conditions for exponential convergence  of the algorithm, while \cite{AMMH:Arxiv:2015} developed a continuous-time algorithm with an exponential convergence rate for a nonsingular and square $A$ and extended the algorithm to the case where $A$ is of full row rank with bounds on the convergence rate.
Furthermore, \cite{LMNM:Automatica:2017} constructed a distributed algorithm and derived the necessary and sufficient conditions on a time-dependent graph for an exponential convergence rate.
Additionally, \cite{SA:TAC:2015} considered distributed computation for a least squares solution to the linear equations that may have no exact solutions, by providing approximate least squares solutions, while \cite{LLAS:Arxiv:2017} dealt with the problem  for the least squares solutions with different graphs and appropriate step-sizes.

Although distributed computation of $Ax=b$ has been studied in the past several years,  the results for distributed computation of general linear matrix equations {
are quite few. }  Note that linear matrix equations  are very important, related to fundamental problems in applied mathematics and computational technology such as the existence of solutions of algebraic equations and stability analysis of  linear systems \cite{BIG:1980,RM:1971}. One of the most famous matrix equations is $AXB=F$ with the matrices $A,\;X,\; B$, and $F$ of appropriate dimensions.   The computation of its solution $X$ plays a fundamental role in  many important application problems such as
the computation of   (generalized) Sylvester equations and generalized inverses of matrices (see \cite{BIG:1980,RM:1971,BK:1980,TW:2013}). It is worthwhile  pointing out that the computation of the special form $AX=F$ (referring to \cite{Woodgate:1996,WC:1996,MHZ:2005}) or a more special form $Ax=b$ as linear algebraic equations with vectors $x$ and $b$ (referring to \cite{SA:TAC:2015,WFM:ACC:2016,SA:ACC:2016,MLM:TAC:2015, AMMH:Arxiv:2015,LLAS:Arxiv:2017,LMM:TAC:2017,LMNM:Automatica:2017}) has also been widely studied for a broad range of applications. 

The objective of this paper is 
to compute   a least squares solution to the well-known matrix equation $AXB=F$ over a multi-agent network in distributed information structures. Considering that the computation of a least squares solution to the linear algebraic equation $Ax=b$ can be related to some optimization problems such as $\min_{x}\|Ax-b\|^2$, we also take a distributed optimization perspective to investigate the solution for this matrix equation over a large-scale network.  Note that  distributed linear matrix equations may have different distributed information structures due to different information structures of $A,$ $B$, and $F$ known by agents.   Based on the column or row sub-blocks  of the matrices $A$, $B$, and $F $ that each agent may know, we  get eight standard matrix-information structures (see Section \ref{Problem} for details), and then provide different substitutional decomposition structures to transform the computation problem to different distributed constrained optimization problems, where each agent only knows local information (instead of the whole data of matrices) and obtains the solution by communicating with its neighbors.
Then we propose distributed continuous-time algorithms and analyze their convergence with help of some control techniques such as stability theory \cite{HHB:TAC:2009} 
and derivative feedback \cite{Antipin:2003}.   In other words, we employ both constrained convex optimization and control ideas to compute a least squares solution to $AXB = F$.
The technical contribution of the paper is summarized as follows.
\begin{itemize}
  \item For a distributed design to solve the linear matrix equation of the form $AXB=F$, we propose eight standard distributed structures, and then construct different decomposition transformations with substitutional variables to reformulate the original computation problem to distributed optimization problems with different constraints (related to consensus intersections or coupled equalities), whose solutions are proved to be  least squares solutions to the original matrix equation. The paper presents a distributed  optimization perspective for investigating distributed computation problems of linear matrix equations.
  \item  Based on the reformulated optimization problems, we design distributed continuous-time algorithms to solve linear matrix equations in the proposed standard structures, respectively, by using modified Lagrangian functions and derivative feedbacks. Because the structures of the problems are different, the proposed algorithms are designed using different techniques in distributed optimization and control although all the algorithms are of primal-dual types. Note that the distributed (continuous-time or discrete-time) algorithms for its very special case $Ax=b$, which were widely investigated in   \cite{SA:TAC:2015,WFM:ACC:2016,SA:ACC:2016,MLM:TAC:2015,CZH:CCC:2017,AMMH:Arxiv:2015,LLAS:Arxiv:2017}, cannot be applied to the computation of the matrix equation.
  \item  For various distributed algorithms in the corresponding structures, we provide rigorous proofs for the correctness and exponential convergence of the algorithms to a   least squares solution based on saddle-point dynamics and stability theory with mild conditions.  Note that some assumptions (such as the existence  of exact solutions or the boundedness of least squares solutions in  \cite{SA:TAC:2015,MLM:TAC:2015,AMMH:Arxiv:2015,CZH:CCC:2017,LMM:TAC:2017,LMNM:Automatica:2017}) for   $Ax=b$ are not required in our paper, and therefore, our results may  also provide a new viewpoint for the  distributed  computation of  $Ax=b$ and its theoretical analysis.
\end{itemize}

The remainder of this paper is organized as follows. Preliminary knowledge is presented in Section \ref{sec:def}, while the problem formulation  of solving a matrix equation with distributed information  and the main result of this paper are given in  Section \ref{Problem}.    Then the  reformulations of the matrix equation in different structures,  distributed algorithms for the reformulated optimization problems, and their exponential convergence   are given in Section \ref{PR}.  Following that, a numerical
simulation is carried out for illustration in Section \ref{sec:num}.  Finally, concluding remarks are provided in Section \ref{conclusion}.

\section{Mathematical Preliminaries}
\label{sec:def}
In this section, we introduce the necessary notations and knowledge related to matrices, graph theory, convex analysis, optimization, and convergence property.

\subsection{Matrices}

Denote $\mathbb R$ as the set of real numbers,
$\mathbb R^n$ as the set of $n$-dimensional real column vectors, $\mathbb R^{n\times m}$ as the set of $n$-by-$m$ real matrices, and $I_n$ as the $n\times n$ identity matrix, respectively.
For $A\in\mathbb R^{m\times n}$, we denote $\mathrm {rank}\,A$ as the rank of  $A$, $A^\mathrm{T}$ as the transpose of $A$, $\mathrm{range} (A)$ as the range of  $A$, $\ker (A)$ as the kernel of $A$, and  ${\rm tr} (A)$ as the trace of $A$.
Write $ 1_{n}$ ($1_{n\times q}$) for the $n$-dimensional column vector ($n\times q$ matrix) with all elements of 1,
$ 0_{n}$ ($0_{n\times q}$) for the $n$-dimensional column  vector ($n\times q$ matrix) with all elements of 0, $A \otimes B$ for the Kronecker product of matrices $A$ and $B$, and ${\rm vec}(A) $ for the vector obtained by stacking the columns of matrix $A$.
Furthermore, denote $\|\cdot\|$  as the Euclidean norm, and $\|\cdot\|_{\rm F}$ as the Frobenius norm of real matrices defined by $\|A\|_{\rm F}=\sqrt{\mathrm{tr}(A^{\rm T}A)}=\sqrt{\sum_{i,j}A_{i,j}^2}$.
Let $\langle \cdot,\cdot \rangle_{\rm F}$ be the Frobenius inner product of real matrices defined by $\langle A_1,A_2 \rangle_{\rm F} = \mathrm{tr}(A_1^{\rm T}A_2)= \sum_{i,j}(A_1)_{i,j}(A_2)_{i,j}$ with $A_1,A_2\in\mathbb R^{m\times n}$, which satisfies $\langle A_1A_2, A_3 \rangle_{\rm F} =\langle A_1, A_3 A_2^{\rm T}\rangle_{\rm F} =\langle A_2, A_1^{\rm T}A_3 \rangle_{\rm F} $ for $A_1\in\mathbb R^{m\times n}$, $A_2\in\mathbb R^{n\times q}$, and $A_3\in\mathbb R^{m\times q}$.
Let $\{m_j\}_{j=1}^n$ and $\{q_j\}_{j=1}^n$ be sequences of $n$ positive integers with $\sum_{j=1}^n m_j=m$ and $\sum_{j=1}^n q_j=q$, and let $A_i\in\mathbb R^{m_i\times q_i}$ for $i\in\{1,\ldots,n\}$. Define  augmented matrices $[A_i ]^{\{m_j\}_{j=1}^n}_{{\rm R}}$ and $[A_i]^{\{q_j\}_{j=1}^n}_{{\rm C}}$ as
\begin{eqnarray}\label{r-ci}
[A_i]^{\{m_j\}_{j=1}^n}_{{\rm R}} &\triangleq & \Big{[} 0_{q_i\times m_1 }, \ldots, 0_{q_i\times m_{i-1} },  A_i^{\rm T},\nonumber\\
&&\ldots, 0_{q_i\times m_{i+1}},\ldots, 0_{q_i\times m_n}
\Big{]} ^{\rm T}\in\mathbb R^{m \times q_i},
\end{eqnarray}
\begin{eqnarray}\label{r-cii}
  [A_i]^{\{q_j\}_{j=1}^n}_{{\rm C}} & \triangleq & \Big{[}  0_{m_i\times  q_1}, \ldots, 0_{m_i\times  q_{i-1}},  A_i, \nonumber\\
  &&\ldots, 0_{m_i\times  q_{i+1}},\ldots, 0_{m_i\times q_n}
\Big{]} \in\mathbb R^{m_i\times q}.
\end{eqnarray}

\subsection{Graph Theory}

An undirected graph $\mathcal G$ is denoted by $\mathcal G(\mathcal V, \mathcal E, A)$, where $\mathcal V=\{1,\ldots, n\}$ is the set of nodes, $\mathcal E\subset\mathcal V \times \mathcal V$ is the set of edges, $ A=[a_{i,j}]\in\mathbb R^{n\times n}$ is the {\em adjacency matrix} such that $a_{i,j}=a_{j,i}>0$ if $\{j,i\}\in\mathcal E$ and $a_{i,j}=0$ otherwise. The {\em Laplacian matrix} is $L_n=D- A$, where $D\in\mathbb R^{n\times n}$ is diagonal with $D_{i,i}=\sum_{j=1}^n a_{i,j}$, $i\in\{1,\ldots,n\}$.
Specifically, if the graph $\mathcal G$ is connected, then
$L_n=L_n^{\rm T}\geq 0$, $\mathrm{rank} \,L_n=n-1$, and $\ker
(L_n)=\{k{1}_n:k\in\mathbb R\}$ \cite{GR2001}.

\subsection{Convex Analysis and Optimization}

A set $\Omega \subseteq \mathbb{R}^p$ is {\em convex} if $\lambda z_1
+(1-\lambda)z_2\in C$ for any $z_1, z_2 \in \Omega$ and $\lambda\in [0,\,1]$.  A function $f: \Omega \to \mathbb{R}$ is said to be {\em convex} (or {\em strictly convex}) if $f(\lambda z_1
+(1-\lambda)z_2) \leq \text{ (or $<$) } \lambda f(z_1) + (1-\lambda)f(z_2)$ for any $z_1, z_2 \in \Omega, z_1 \neq z_2$ and $\lambda\in (0,\,1)$.
Sometimes, a convex optimization problem can be written as
$\min_{z\in\Omega} f(z),$
where $\Omega \subseteq \mathbb{R}^p$ is a convex set and $f:\mathbb R^p\rightarrow\mathbb R$ is a convex function.

\subsection{Convergence Property}

Consider a dynamical system
\begin{eqnarray}\label{DS}
  \dot x(t) = \phi (x(t)),\quad x(0)=x_0,\quad t\geq 0,
\end{eqnarray}
where $\phi:\mathbb R^q \rightarrow\mathbb R^q$ is Lipschitz continuous. Given a trajectory $x:[0,\infty)\rightarrow \mathbb R^q$ of \eqref{DS},  $y$ is a positive limit point of $x(\cdot)$ if there is a positive increasing divergent sequence $\{t_i\}_{i=1}^\infty\subset \mathbb R$ such that $y=\lim_{i\rightarrow\infty}x(t_i)$, and a positive limit set of $x(\cdot )$ is the set of all positive limit points of $x(\cdot)$.
A set $\mathcal D$ is said to be {\em positive invariant} with respect to (\ref{DS}) if  $x(t)\in\mathcal D$ for all $t\geq 0$ and every $x_0\in\mathcal D$.

Denote $\mathcal B_{\epsilon}(x),x\in\mathbb R^n$ with a constant $\epsilon>0$ as the open ball {\em centered} at $x$ with {\em radius} $\epsilon$.
Let $\mathcal D\subset\mathbb R^q$ be a positive invariant set with respect to \eqref{DS} and $z\in\mathcal D$ be an equilibrium of \eqref{DS}.
$z$ is {\em Lyapunov stable} if, for every $\epsilon>0$, there exists $\delta=\delta(\epsilon) >0$ such that, for every initial condition $x_0\in {\mathcal B_{\delta}(z)\bigcap\mathcal D}$,     the solution $x(t)$ of \eqref{DS} stays in $\mathcal B_{\epsilon}(z)$ for all $t\geq 0$.

The following Lemmas are needed in the analysis of this paper.

\begin{lemma}\cite[Theorem 3.1]{HHB:TAC:2009}\label{semistability}
Let $\mathcal D$ be a compact, positive invariant set with respect to system (\ref{DS}), $V:\mathbb R^q\rightarrow\mathbb R$ be a continuously differentiable function, and $ x(\cdot) \in\mathbb R^{q}$ be a solution of (\ref{DS}) with $ x(0) = x_0\in\mathcal D$. Assume
$\dot V(x)\leq 0,\quad \forall x\in\mathcal D,$
and define $\mathcal Z=\{x\in\mathcal D:\dot V(x)=0\}.$
If every point in the largest invariant subset $\mathcal M$ of {$\overline{\mathcal Z}\bigcap\mathcal D$} is Lyapunov stable, where $\overline{\mathcal{Z}}$ is the closure of $\mathcal{Z}\subset\mathbb{R}^{n}$, then (\ref{DS}) converges to one of its Lyapunov stable equilibria for any $x_0\in\mathcal D$.
\end{lemma}

\begin{lemma}\label{exponential}
Suppose  $\phi (x) =  Mx+b$ with $M\in\mathbb R^{q\times q}$, $b\in\mathbb R^q$, and $\mathcal D=\mathbb R^q$. The following statements are equivalent.
\begin{itemize}
  \item [(i)] System (\ref{DS}) converges to an equilibrium exponentially for any initial condition.
  \item [(ii)] System (\ref{DS}) converges to an equilibrium  for  any initial condition.
\end{itemize}
\end{lemma}
\begin{IEEEproof}
It is trivial that (i)$\Rightarrow$(ii) and the proof is omitted.

Suppose (ii) holds. Let $x^*$ be an equilibrium of (\ref{DS}) and define $y = x-x^*$.   System (\ref{DS}) is equivalent to
\begin{eqnarray}\label{equiv1}
\dot y(t) = My(t),\quad y(0)=y_0\in \mathbb R^q,\quad  t\geq 0.
\end{eqnarray}
To show (ii)$\Rightarrow$(i), we show that $y(t)$  converges to an equilibrium exponentially for  any initial condition.

It follows from statement (ii) and Definition 11.8.1 of \cite[p. 727]{Bernstein:2009} that $M$ is semistable (that is, its eigenvalues lie on the open left half complex plane, except for a few semi-simple zero eigenvalues). Hence, there exists an invertible matrix $P\in\mathbb R^{q\times q}$ such that $PMP^{-1}=\begin{bmatrix}
&D &0_{r \times (q-r)}\\
&0_{(q-r)\times r} & 0_{(q-r)\times (q-r)}
\end{bmatrix}$, where $D\in\mathbb R^{r\times r}$ is Hurwitz and $r=\mathrm{rank}\,M\leq q$.
Define $z=\begin{bmatrix}
z_1 \\
 z_2
\end{bmatrix}=Py\in\mathbb R^q$ such that $z_1\in\mathbb R^{r}$, and $z_2\in\mathbb R^{q-r}$. It follows from \eqref{equiv1} that
\begin{eqnarray*}
\dot z_1(t) = Dz_1(t),\, \dot z_2(t)=0_{q-r},\, z_1(0)=z_{1,0},\, z_2(0)=z_{2,0},\,  t\geq 0.
\end{eqnarray*}
Hence, $\Bigg{\|}z(t)-\begin{bmatrix}
0_r \\
 z_2(0)
\end{bmatrix}\Bigg{\|} = \|z_1(t)\|=\|\mathrm{e}^{Dt}z_1(0)\|$. Recall that $D$ is Hurwitz. The trajectory  $z(t)$ converges to $\begin{bmatrix}
0_r \\
 z_2(0)
\end{bmatrix}$ exponentially and, equivalently, $y(t)$ converges to $P^{-1}\begin{bmatrix}
0_r \\
 z_2(0)
\end{bmatrix}$ exponentially.
\end{IEEEproof}

\section{Problem Description and Main Result}
\label{Problem}
In this paper, we consider the distributed computation of a least squares solution to the well-known matrix equation in the following form
\begin{eqnarray}\label{AXBeF}
AXB=F,
\end{eqnarray}
where  $A\in\mathbb R^{m\times r}$,  $B\in\mathbb R^{p\times q}$, and $F\in\mathbb R^{m\times q}$ are known matrices, and $X\in\mathbb R^{r\times p}$ is an unknown matrix to be solved.
Note that equation \eqref{AXBeF} may not have a solution $X$. However, it always has a least squares solution, which is defined as follows.

\begin{definition}
A {\em least squares solution} to \eqref{AXBeF} is a solution of the optimization problem $
\min_{X}\|AXB-F\|_{\rm F}^2.
$
\end{definition}

Obviously, if \eqref{AXBeF} has a solution, then a least squares solution is also an  exact solution.
The following result is well known (see \cite{DLD:2008,YonggeTian:2009}).
\begin{lemma}~
\begin{itemize}
  \item [1.] Equation \eqref{AXBeF} has an exact solution if and only if
        $$\mathrm{range}(F)\subset\mathrm{range}(A)\text{ and }\mathrm{range}(F^{\rm T})\subset\mathrm{range}(B^{\rm T}).$$
  \item [2.] $X^*\in\mathbb R^{r\times p}$ is a least squares solution if and only if
        \begin{eqnarray}\label{opt-cond-c}
        0_{r\times p}= \frac{\partial \|AXB-F\|_{\rm F}^2}{\partial X}\Big{|}_{X=X^*}=A^{\rm T}(AX^*B-F)B^{\rm T}.
        \end{eqnarray}
  \item [3.] If $A$ is full column-rank and $B$ is full row-rank, $X^*=(A^{\rm T}A)^{-1}A^{\rm T}FB^{\rm T}(BB^{\rm T})^{-1}$ is the unique least squares solution.
\end{itemize}
\end{lemma}
Note that \eqref{AXBeF} is one of the most famous matrix equations in matrix theory and applications  (see \cite{BIG:1980,RM:1971}), related to the computation of many important problems such as (generalized) Sylvester equations and generalized inverses of matrices (see \cite{BIG:1980,RM:1971,BK:1980,DLD:2008,TW:2013}). Because solving \eqref{AXBeF} is one of the key problems of matrix computation, many techniques have been proposed and various  centralized algorithms have been developed to solve problem   \eqref{AXBeF} (see \cite{DLD:2008,PHZ:2005,TW:2013,HPZ:2006,ZLGZ:2011}).   One significant method is a gradient-based approach from the optimization viewpoint (see Theorem 2 of \cite{DLD:2008}). Because many matrix equations in engineering and science fields have large scales, the distributed computation of  \eqref{AXBeF} is very necessary.
However, {
very few results} have been obtained for the distributed computation of  \eqref{AXBeF} due to its complicated structures when each agent only knows some sub-blocks of  (large-size) matrices $A$, $B$, and $F$.

On the other hand, distributed computation of linear algebraic equations in the form of $Ax=b$ with vectors $x$ and $b$ has been widely studied in recent years and some significant results have been obtained in \cite{SA:TAC:2015,WFM:ACC:2016,SA:ACC:2016,MLM:TAC:2015, AMMH:Arxiv:2015,LLAS:Arxiv:2017,CZH:CCC:2017}.   To solve  \eqref{AXBeF}, an immediate idea is to vectorize it as follows:
$$
\mathrm{vec}(AXB) = (B^{\rm T}\otimes A)\mathrm{vec}(X)=\mathrm{vec}(F),
$$
and try the existing linear algebraic equation results here.   Although this idea may work in centralized situations, it may totally spoil the original distributed information structure because
the local knowledge about some sub-blocks of $A$ and $B$ of each agent may be mixed up and multiplied due to the Kronecker product.
Hence, we have to develop new methods to   solve the matrix equation  \eqref{AXBeF} in a distributed way.

In this paper, we consider the distributed computation of a least squares solution to \eqref{AXBeF} over a multi-agent network described by an undirected  graph $\mathcal G$, where matrices $A$, $B$, and $F$ are composed of $n$ row-block or column-block matrices, known by $n$ agents.

In this complicated problem, there are different distributed information structures of matrices $A$, $B$, and $F$. To distinguish the row-blocks or column-blocks of a matrix, we use subscript ``$vi$" to denote its $i$th row-block and subscript ``$li$" to denote its $i$th column-block in the sequel.

For different information structures of matrices $A$, $B$, and $F$, we can classify the distributed computation problem of equation \eqref{AXBeF} in the following {\bf eight} standard structures:
\begin{itemize}
\item {\bf Row-Column-Column (RCC) Structure}:  Consider \eqref{AXBeF} with
       \begin{eqnarray}\label{RCC}
            &&A = \begin{bmatrix} A_{v1}\\
            \vdots\\
            A_{vn}
            \end{bmatrix}\in\mathbb R^{m\times r}, \quad B = \begin{bmatrix} B_{l1},
            \ldots,
            B_{ln}
            \end{bmatrix}\in\mathbb R^{p\times q}, \nonumber\\
            && F = \begin{bmatrix} F_{l1},
            \ldots,
            F_{ln}
            \end{bmatrix}\in\mathbb R^{m\times q},
        \end{eqnarray}
where $A_{vi}\in\mathbb R^{m_i\times r}$,  $B_{li}\in \mathbb R^{p\times q_i}$, $F_{li}\in \mathbb R^{m\times q_i}$, $\sum_{i=1}^n m_i = m$,  $\sum_{i=1}^n q_i = q$,
and the sub-blocks of $A$, $B$, and $F$ are distributed among the agents of network $\mathcal G$.

In this structure, agent $i$ only knows  $A_{vi}$, $B_{li}$, and $F_{li}$ for $i\in\{1,\ldots,n\}$.
By communicating with its neighbors, every agent $i$ aims to obtain a least squares solution to \eqref{AXBeF}.

Obviously, if $X$ and $F$ are row vectors with $A=1$,
the matrix equation   \eqref{AXBeF} with \eqref{RCC} becomes a linear algebraic equation $B^{\rm T}X^{\rm T} = F^{\rm T}$, where each agent knows a  row sub-block of $B^{\rm T}$ and the vector $F^{\rm T}$, which was investigated in \cite{MLM:TAC:2015,SA:TAC:2015,WFM:ACC:2016,SA:ACC:2016,LMM:TAC:2017,LMNM:Automatica:2017} and references therein. However, the sub-blocks of matrices $A$, $B$, and $X$ are coupled in the original equation \eqref{AXBeF} with \eqref{RCC}, and hence, new techniques and ideas are needed for its distributed algorithm design.

\item {\bf Row-Row-Row (RRR) Structure}:  Consider \eqref{AXBeF} with
\begin{eqnarray}\label{RRR}
            && A = \begin{bmatrix} A_{v1}\\
            \vdots\\
            A_{vn}
            \end{bmatrix}\in\mathbb R^{m\times r},
            \quad B = \begin{bmatrix} B_{v1}\\
            \vdots\\
            B_{vn}
            \end{bmatrix}\in\mathbb R^{p\times q}, \nonumber\\
            && F = \begin{bmatrix} F_{v1}\\
            \vdots\\
            F_{vn}
            \end{bmatrix}\in\mathbb R^{m\times q},
        \end{eqnarray}
with
$ X = [X_{l1},\ldots,X_{ln}]\in\mathbb R^{r\times p},$
where $A_{vi}\in\mathbb R^{m_i\times r}$, $X_{li}\in\mathbb R^{r\times p_i}$, $B_{vi}\in \mathbb R^{p_i\times q}$, $F_{vi}\in \mathbb R^{m_i\times q}$, $\sum_{i=1}^n m_i = m$, and $\sum_{i=1}^n p_i = p$.  Similarly, agent $i$ in the $n$-agent network $\mathcal G$ only knows $A_{vi}$, $B_{vi}$, and $F_{vi}$ and cooperates with its neighbors to compute $X_{li}$.

Clearly,
if $X$ and $ F$ are  row vectors with $A=1$, this problem becomes that discussed in \cite{CZH:CCC:2017}.

\item {\bf Column-Column-Row (CCR) Structure}:
Consider \eqref{AXBeF} with
    \begin{eqnarray}\label{CCR}
            && A = \begin{bmatrix} A_{l1}, \ldots, A_{ln} \end{bmatrix}\in\mathbb R^{m\times r},
            \quad  F = \begin{bmatrix} F_{v1}\\
            \vdots\\
            F_{vn}
            \end{bmatrix}\in\mathbb R^{m\times q},
            \nonumber\\
            &&B = \begin{bmatrix} B_{l1},
            \ldots,
            B_{ln}
            \end{bmatrix}\in\mathbb R^{p\times q},
      \end{eqnarray}
where $A_{li}\in\mathbb R^{m\times r_i}$,  $B_{li}\in \mathbb R^{p\times q_i}$, $F_{vi}\in \mathbb R^{m_i\times q}$, $\sum_{i=1}^n r_i = r$, $\sum_{i=1}^n m_i = m$, and $\sum_{i=1}^n q_i = q.$
We use an $n$-agent network $\mathcal G$ to find $X$, where agent $i$ knows $A_{li}$, $B_{li}$, and $F_{vi}$ and estimates $X$ by cooperating with its neighbors to reach a consensus to a least squares solution to matrix equation  \eqref{AXBeF} with \eqref{CCR}.

\item{\bf Column-Row-Row (CRR) Structure}:  Consider \eqref{AXBeF} with
       \begin{eqnarray}\label{CRR}
            && A = \begin{bmatrix} A_{l1}, \ldots, A_{ln} \end{bmatrix}\in\mathbb R^{m\times r}, \quad B = \begin{bmatrix} B_{v1}\\
            \vdots\\
            B_{vn}
            \end{bmatrix}\in\mathbb R^{p\times q}, \nonumber\\
            && F = \begin{bmatrix} F_{v1}\\
            \vdots\\
            F_{vn}
            \end{bmatrix}\in\mathbb R^{m\times q},
        \end{eqnarray}
where $A_{li}\in\mathbb R^{m\times r_i}$, $X_{li}\in\mathbb R^{r\times p_i}$, $X = [X_{l1},\ldots,X_{ln}]\in\mathbb R^{r\times p} $, $B_{vi}\in \mathbb R^{p_i\times q}$, $F_{vi}\in \mathbb R^{m_i\times q}$, $\sum_{i=1}^n r_i = r$, $\sum_{i=1}^n m_i = m$, and $\sum_{i=1}^n p_i = p$.
We use an $n$-agent system to find $X$, where agent $i$ knows $A_{li}$, $B_{vi}$, and $F_{vi}$ and cooperates with its neighbors to compute $X_{li}$, which composes a least squares solution to matrix equation \eqref{AXBeF} and \eqref{CRR}.

In this case, if  $X$ and $ F$ are  row vectors with $A=1$, \eqref{AXBeF} and \eqref{CRR}  becomes the problem investigated in \cite{CZH:CCC:2017}.

\item {\bf Row-Column-Row (RCR) Structure}: Consider  \eqref{AXBeF} with
    \begin{eqnarray}\label{RCR}
           && A = \begin{bmatrix} A_{v1}\\
            \vdots\\
            A_{vn}
            \end{bmatrix}\in\mathbb R^{m\times r}, \, B = \begin{bmatrix} B_{l1},
            \ldots,
            B_{ln}
            \end{bmatrix}\in\mathbb R^{p\times q}, \nonumber\\
            && F = \begin{bmatrix} F_{v1}\\
            \vdots\\
            F_{vn}
            \end{bmatrix}\in\mathbb R^{m\times q}.
    \end{eqnarray}
Clearly, this structure is equivalent to RCC structure by the transposes of matrices.

\item {\bf Column-Column-Column (CCC) Structure}: Consider  \eqref{AXBeF} with
      \begin{align}\label{CCC}
            A &= \begin{bmatrix} A_{l1}, \ldots, A_{ln} \end{bmatrix}\in\mathbb R^{m\times r},
          \, B = \begin{bmatrix} B_{l1},
            \ldots,
            B_{ln}
            \end{bmatrix}\in\mathbb R^{p\times q},\nonumber\\
            F &= \begin{bmatrix} F_{l1},
            \ldots,
            F_{ln}
            \end{bmatrix}\in\mathbb R^{m\times q}.
      \end{align}
It is equivalent to  RRR structure by the transposes of matrices.

\item {\bf Row-Row-Column (RRC) Structure}:  Consider \eqref{AXBeF} with
        \begin{eqnarray}\label{RRC}
          &&  A = \begin{bmatrix} A_{v1}\\
            \vdots\\
            A_{vn}
            \end{bmatrix}\in\mathbb R^{m\times r},
            \quad B = \begin{bmatrix} B_{v1}\\
            \vdots\\
            B_{vn}
            \end{bmatrix}\in\mathbb R^{p\times q}, \nonumber\\
            && F = \begin{bmatrix} F_{l1},
            \ldots,
            F_{ln}
            \end{bmatrix}\in\mathbb R^{m\times q}.
        \end{eqnarray}
It is equivalent to  CCR structure by the transposes of matrices.

\item {\bf Column-Row-Column (CRC) Structure}: Consider \eqref{AXBeF} with
        \begin{eqnarray}\label{CRC}
           && A = \begin{bmatrix} A_{l1}, \ldots, A_{ln} \end{bmatrix}\in\mathbb R^{m\times r},
            \quad B = \begin{bmatrix} B_{v1}\\
            \vdots\\
            B_{vn}
            \end{bmatrix}\in\mathbb R^{p\times q}, \nonumber\\
            && F = \begin{bmatrix} F_{l1},
            \ldots,
            F_{ln}
            \end{bmatrix}\in\mathbb R^{m\times q}.
        \end{eqnarray}
It is equivalent to  CRR structure by the transposes of matrices.
\end{itemize}

{\begin{remark}
In  these formulations, the rows and columns of matrices may have different physical interpretations. Take $A$ for example. If $A$ is decomposed of row blocks, each row defines a local  linear space and matrix $A$ defines the intersection of the local linear spaces as in \cite{MLM:TAC:2015,SA:TAC:2015,WFM:ACC:2016,SA:ACC:2016,LMM:TAC:2017,LMNM:Automatica:2017}. However, if $A$ is decomposed of column blocks, each block contains partial information of the coupling/monotropic constraint information, for example in resource allocation problems as discussed in \cite{YHL:SCL:2015}. Due to  different  structures,  different combinations of the consensus design and the auxiliary decomposition design are adopted.
\hfill$\Diamond$
\end{remark}}

The main result of this paper can, in fact, be written as

\begin{theorem}\label{themain}
{\bf A least squares solution to \eqref{AXBeF} in the eight standard structures can be obtained using distributed algorithms with exponential convergence rates if the undirected graph $\mathcal G$ is connected}.
\end{theorem}

Clearly, because   RCR, CCC, RRC, and CRC structures are the transpose of RCC, RRR, CCR, and CRR structures,
the eight different structures are basically four standard structures in the distributed computation design.  Therefore, we only need to study \eqref{AXBeF} with the four standard structures, \eqref{RCC}-\eqref{CRR}, in the sequel.

\section{Reformulation, Algorithm, and Exponential Convergence}\label{PR}
In this section,  we first reformulate the matrix computation problem in four different structures as solvable distributed optimization problems with different substitutional decompositions.
Then we propose distributed continuous-time algorithms for the four standard structures using a derivative feedback idea and the saddle-point dynamics. Finally, we give the exponential convergence proof of our algorithms with help of the stability theory and the Lyapunov method.

\subsection{Row-Column-Column Structure}
To handle the couplings between the sub-blocks of matrices $A$, $B$, and $F$ in the equation \eqref{AXBeF} with \eqref{RCC},
we introduce a substitutional variable $Y$ to make \eqref{AXBeF} and \eqref{RCC} equivalent to $Y = AX$ and $YB_{li} = F_{li}$ for $i\in\{1,\ldots,n\}$.
Let $X_i\in\mathbb R^{r\times p}$ and $Y_i\in\mathbb R^{m\times p}$  be the estimates of $X$ and $Y$ of agent $i\in\{1,\ldots,n\}$,  respectively.  We propose a {\bf full-consensus substitutional decomposition} method by requiring both $X_i$ and $Y_i$ to achieve consensus, namely, we rewrite the equation \eqref{AXBeF} with \eqref{RCC} as
\begin{eqnarray}
Y_iB_{li} &=& F_{li}, \quad Y_i=Y_j,\quad i,j\in\{1,\ldots,n\},
\label{RCC-DP1}\\
AX_i  &=& Y_i,\quad X_i=X_j. \label{RCC-DP2}
\end{eqnarray}

Clearly, \eqref{RCC-DP2} is not in a distributed form because all the sub-blocks of $A$ need to be known. To decompose  \eqref{RCC-DP2}, define $Y_i \triangleq \begin{bmatrix}
Y_{i}^{v1}\\
\vdots\\
Y_{i}^{vn}
\end{bmatrix}$, where  $Y_{i}^{vj}\in\mathbb R^{m_j\times p}$ for all $i,j\in\{1,\ldots,n\}$.  Due to $Y_i=Y_j$ for all $i,j\in\{1,\ldots,n\}$ in \eqref{RCC-DP1},
\eqref{RCC-DP2} is equivalent to
\begin{eqnarray}
A_{vi} X_i = Y_{i}^{vi},\quad X_i=X_j,\quad i,j\in\{1,\ldots,n\}.\label{RCC-DP3}
\end{eqnarray}

Hence, the matrix equation \eqref{AXBeF} with \eqref{RCC} is equivalent to the linear matrix equations \eqref{RCC-DP1} and \eqref{RCC-DP3}.   Define  extended matrices $X_{\rm E} = [X_1^{\rm T},\cdots, X_n^{\rm T}]^{\rm T}\in\mathbb R^{nr\times p}$ and $Y_{\rm E}=[Y_1^{\rm T},\cdots, Y_n^{\rm T}]^{\rm T}\in\mathbb R^{nm\times p}$.  Based on \eqref{RCC-DP1} and \eqref{RCC-DP3}, we {\bf reformulate} the distributed computation of \eqref{AXBeF} with RCC structure as   the following distributed optimization problem
\begin{subequations}\label{opt-rcc}
\begin{align}
\min_{X_{\rm  E},Y_{\rm E}}&\quad {\sum_{i=1}^n\|Y_iB_{li} - F_{li}\|_{\rm F}^2}, \label{opt-rcc1}\\
\text{s. t.}&\quad X_i=X_j,\, Y_i =Y_j, \, A_{vi}X_i  = Y_i^{vi},\,i,j\in\{1,\ldots,n\},\label{opt-rcc2}
\end{align}
\end{subequations}
where agent $i$ knows $A_{vi}$, $B_{li}$, $F_{li}$, and estimates the solution $X_i$ and $Y_i$ with only local information.

\begin{remark}
Problem \eqref{opt-rcc} is a standard distributed optimization problem, which contains local constraints $A_{vi}X_i  = Y_i^{vi}$ and consensus constraints $X_i=X_j$ and $Y_i =Y_j$.
\hfill$\Diamond$
\end{remark}

The following proposition reveals the relationship between \eqref{AXBeF} and problem \eqref{opt-rcc}.

\begin{proposition}\label{rcc-se}
Suppose that the undirected graph $\mathcal G$ is connected.   $X^*\in\mathbb R^{r\times p}$ is a least squares solution to matrix equation  \eqref{AXBeF} if and only if
there exists $Y^*=AX^*\in\mathbb R^{m\times p}$ such that $(X_{\rm E}^*,Y_{\rm E}^*)=(1_n\otimes X^*,1_n\otimes Y^*)$  is a solution to problem \eqref{opt-rcc}.
\end{proposition}

The proof can be found in Appendix \ref{app-rcc-se}.

In this structure, we focus on problem \eqref{opt-rcc}, and propose a distributed {\bf algorithm} of agent $i$ as
\begin{subequations}\label{algorithm2}
\begin{align}
\dot X_i(t) &= -A_{vi}^{\rm T}(A_{vi}X_i(t)-Y_i^{vi}(t))-A_{vi}^{\rm T}\Lambda^3_i(t)\nonumber\\
&\quad- \sum_{j=1}^n a_{i,j}(\Lambda^1_i(t)-\Lambda^1_j(t))\nonumber\\
&\quad-\sum_{j=1}^n a_{i,j}(X_i(t)-X_j(t)),\, X_i(0)=X_{i0}\in\mathbb R^{r\times p},\label{alm2-x}\\
\dot Y_i(t) &= -(Y_i(t)B_{li} - F_{li})B_{li}^{\rm T}+{[I_{m_i}]_{\rm R} \Lambda^3_i}(t) \nonumber\\
&\quad+[I_{m_i}]_{\rm R}(A_{vi}X_i(t)-Y_i^{vi}(t))-\sum_{j=1}^n a_{i,j}(Y_i(t)-Y_j(t))\nonumber\\
&\quad -\sum_{j=1}^n a_{i,j}(\Lambda^2_i(t)-\Lambda^2_j(t)),\, Y_i(0)=Y_{i0}\in\mathbb R^{m\times p},\label{alm2-y}\\
\dot \Lambda^1_i(t) &= \sum_{j=1}^n a_{i,j}(X_i(t)-X_j(t)),\quad \Lambda^1_i(0)=\Lambda^1_{i0}\in\mathbb R^{r\times p},\\
\dot \Lambda^2_i(t) &= \sum_{j=1}^n a_{i,j}(Y_i(t)-Y_j(t)),\quad \Lambda^2_i(0)=\Lambda^2_{i0}\in\mathbb R^{m\times p},\\
\dot \Lambda^3_i(t) &= A_{vi}X_i(t)-Y_i^{vi}(t),\quad \Lambda^3_i(0)=\Lambda^3_{i0}\in\mathbb R^{m_i\times p},
\end{align}
\end{subequations}
where $i\in\{1,\ldots,n\}$, $t\geq 0$,
$X_i(t)$ and $Y_i(t)$ are the estimates of solutions to problem \eqref{opt-rcc} by agent $i$ at time $t$, $\Lambda^1_i(t)$, $\Lambda^2_i(t)$, and $\Lambda^3_i(t)$ are the estimates of Lagrangian multipliers for the constraints in \eqref{opt-rcc2} by agent $i$ at time $t$, and $[I_{m_i}]_{\rm R}$ denotes $[I_{m_i}]^{\{m_j\}_{j=1}^n}_{{\rm R}}$, as defined in \eqref{r-ci}.

\begin{remark}
Algorithm \eqref{algorithm2} is a primal-dual algorithm, whose primal variables are $X_i$ and $Y_i$ and dual variables are $\Lambda^1_i$, $\Lambda^2_i$, and $\Lambda^3_i$. Though  substitutional variables  are used in \eqref{algorithm2} for the distributed computation of \eqref{AXBeF} and \eqref{RCC}, algorithm \eqref{algorithm2} is a fully distributed algorithm. {Different from the classic (Arrow-Hurwicz-Uzawa type) primal-dual algorithm in \cite{LLAS:Arxiv:2017}, the consensus design ($-\sum_{j=1}^n a_{i,j}(X_i(t)-X_j(t))$ and $ -\sum_{j=1}^n a_{i,j}(Y_i(t)-Y_j(t))$ in \eqref{alm2-x} and \eqref{alm2-y})
and the damping design ($-A_{vi}^{\rm T}(A_{vi}X_i(t)-Y_i^{vi}(t))$ in \eqref{alm2-x}) are used in the algorithm.}
\hfill$\Diamond$
\end{remark}

\begin{remark}
Let  $\Lambda^1 = \begin{bmatrix}\Lambda^1_1\\
\vdots\\
\Lambda^1_n
\end{bmatrix}\in\mathbb R^{nr\times p}$, $\Lambda^2 = \begin{bmatrix}\Lambda^2_1\\
\vdots\\
\Lambda^2_n
\end{bmatrix}\in\mathbb R^{nm\times p}$, and $\Lambda^3 = \begin{bmatrix}\Lambda^3_1\\
\vdots\\
\Lambda^3_n
\end{bmatrix}\in\mathbb R^{m\times p}$, where $\Lambda^1_i\in\mathbb R^{r\times p}$, $\Lambda^2_i\in\mathbb R^{m\times p}$, and $\Lambda^3_i\in\mathbb R^{m_i\times p}$. Algorithm \eqref{algorithm2} can be viewed as the saddle-point dynamics of the modified  Lagrangian function $L(X_{\mathrm E},Y_{\mathrm E},\Lambda^1,\Lambda^2,\Lambda^3)$ given by
\begin{align*}
L =& \frac{1}{2}\sum_{i=1}^n\|Y_iB_{li} - F_{li}\|_{\rm F}^2
+\sum_{i=1}^n\sum_{j=1}^n \Big\langle \Lambda^1_i, a_{i,j}(X_i-X_j)\Big\rangle_{\rm F}\\
&+\sum_{i=1}^n\sum_{j=1}^n\Big\langle \Lambda^2_i,  a_{i,j}(Y_i-Y_j)\Big\rangle_{\rm F} +\sum_{i=1}^n\langle \Lambda^3_i, A_{vi}X_i-Y_i^{vi} \rangle_{\rm F}\\
&+\frac{1}{2}\sum_{i=1}^n\sum_{j=1}^n\Big\langle X_i,  a_{i,j}(X_i-X_j)\Big\rangle_{\rm F}
 +\frac{1}{2}\sum_{i=1}^n\|A_{vi}X_i-Y_i^{vi}\|^2_{\rm F}\\
 &+\frac{1}{2}\sum_{i=1}^n\sum_{j=1}^n \Big\langle Y_i, a_{i,j}(Y_i-Y_j)\Big\rangle_{\rm F},
\end{align*}
where $a_{i,j}$ is the $(i,j)$th element of the adjacency matrix of graph $\mathcal G$, $\Lambda^1$, $\Lambda^2$, and $\Lambda^3$  are the Lagrangian matrix multipliers, that is,
$\dot X_i = -\nabla_{X_i} L,$ $\dot Y_i = -\nabla_{Y_i} L,$ $\dot \Lambda^1_i = \nabla_{\Lambda^1_i} L,$ $\dot \Lambda^2_i =  \nabla_{\Lambda^2_i} L,$ and $\dot \Lambda^3_i =  \nabla_{\Lambda^3_i} L$ for $i\in\{1,\ldots,n\}$.  In this modified Lagrangian function, $\frac{1}{2}\sum_{i=1}^n\|A_{vi}X_i-Y_i^{vi}\|^2_{\rm F}$ is the augmented term, and $\frac{1}{2}\sum_{i=1}^n\langle X_i, \sum_{j=1}^n a_{i,j}(X_i-X_j)\rangle_{\rm F}$ and $\frac{1}{2}\sum_{i=1}^n\langle Y_i, \sum_{j=1}^n a_{i,j}(Y_i-Y_j)\rangle_{\rm F}$ are the (weighted) Laplacian regularization.
\hfill$\Diamond$
\end{remark}

\begin{remark}
Recall that a standard centralized algorithm (see \cite{DLD:2008}) is $X(k)=X(k-1)+\mu A^{\rm T}[F-AX(k-1)B]B^{\rm T},$
where  $\mu>0$ is an appropriate real number and $A^{\rm T}[F-AX(k-1)B]B^{\rm T}$ is the negative gradient  of $\frac{1}{2}\|AXB-F\|_{\rm F}^2$. Both the centralized algorithm and algorithm \eqref{algorithm2} use the gradient dynamics in the design. However, in contrast to the centralized algorithm, algorithm \eqref{algorithm2} uses auxiliary variables and equality linear constraints to deal with the unavailability of matrices information. The auxiliary variables may also be considered as distributed observers
and filters of the unavailable matrix information from the control viewpoint.
\hfill$\Diamond$
\end{remark}

The following result reveals the relationship between equilibria of algorithm \eqref{algorithm2} and solutions to problem \eqref{opt-rcc}, which is an immediate conclusion of the KKT optimality condition (Theorem 3.25 of \cite{Ruszczynski:2006}), so its proof is omitted.

\begin{lemma}\label{rcc-eq}
For a connected undirected graph $\mathcal G$, $(X_{\rm E}^*,Y_{\rm E}^*)\in\mathbb R^{nr\times p}\times \mathbb R^{nm\times p}$ is a solution to problem \eqref{opt-rcc} if and only if there exist matrices $\Lambda^{1*}\in\mathbb R^{nr\times p}$, $\Lambda^{2*}\in\mathbb R^{nm\times p}$, and $\Lambda^{3*}\in\mathbb R^{m\times p}$ such that $(X_{\rm E}^*,Y_{\rm E}^*,\Lambda^{1*},\Lambda^{2*},\Lambda^{3*})$  is an equilibrium of \eqref{algorithm2}.
\end{lemma}

It is time to show the {\bf exponential convergence} of algorithm \eqref{algorithm2}.

\begin{proposition}\label{THM_RCC}
If the undirected graph $\mathcal G$ is connected,  then
\begin{enumerate}
  \item every equilibrium  of algorithm \eqref{algorithm2} is Lyapunov stable and its trajectory is bounded for any initial condition;

  \item {the trajectory of algorithm \eqref{algorithm2} is exponentially convergent and} $X_i(t)$ converges to a least squares solution to \eqref{AXBeF} exponentially for all $i\in\{1,\ldots,n\}$.
\end{enumerate}
\end{proposition}

The proof can be found in Appendix \ref{app-THM_RCC}, which  shows that algorithm \eqref{algorithm2} is globally {convergent}.  In fact, if there are multiple solutions, the solution obtained by algorithm \eqref{algorithm2} depends on the selected initial condition.

\begin{remark}
In comparison with many previous results on the distributed computation of linear algebraic equation $Ax=b$, we do not need the boundedness assumption for least squares solutions given in \cite{SA:TAC:2015} or the existence assumption of  exact solutions given in \cite{MLM:TAC:2015, AMMH:Arxiv:2015,CZH:CCC:2017,LMM:TAC:2017,LMNM:Automatica:2017}. By virtue of  the saddle-point dynamics, \cite{WZMC:arx:2017} proposed a distributed algorithm  to achieve least squares solutions with an exponential convergence rate, without any boundedness assumption or  choices of time-varying small step sizes.
\hfill $\Diamond$
\end{remark}

\subsection{Row-Row-Row Structure}

The decomposition method used in RCC structure cannot convert RRR structure to a solvable optimization problem.   To deal with \eqref{AXBeF} and \eqref{RRR}, we take substitutional variables $Y_i\in\mathbb R^{r\times q}$ such that $Y_i = XB$ for all $i\in\{1,\ldots,n\}$, and then we propose another method called {\bf $Y$-consensus substitutional decomposition} because we need the consensus of $Y_i$ as follows:
\begin{eqnarray}
A_{vi}Y_i &=&F_{vi}, \quad Y_i=Y_j,\quad i,j\in\{1,\ldots,n\},    \label{RRR-DP1}\\
\frac{1}{n}\sum_{i=1}^{n} Y_i &=& [X_{l1},\ldots,X_{ln}]\begin{bmatrix} B_{v1}\\
\vdots\\
B_{vn}
\end{bmatrix} = \sum_{i=1}^{n}X_{li}B_{vi} .   \label{RRR-DP2}
\end{eqnarray}

In this way, agent $i$ computes $X_{li}$ and $Y_i$ based on only local information.
To decompose \eqref{RRR-DP2}, we add new variables $Z_i\in\mathbb R^{r\times q}$ such that
\begin{eqnarray}\label{rrr-extend}
\frac{1}{n}Y_i-X_{li}B_{vi} +\sum_{j=1}^n a_{i,j}(Z_i-Z_j) = 0_{r \times q},
\end{eqnarray}
where $a_{i,j}$ is the $(i,j)$th element of the adjacency matrix of the connected graph for the agents. Clearly, \eqref{rrr-extend} implies \eqref{RRR-DP2}. Conversely, if \eqref{RRR-DP2} holds, there exists $Z_i\in\mathbb R^{r\times q}$ such that \eqref{rrr-extend} holds due to the fundamental theorem of linear algebra \cite{Strang:1993} (whose proof is similar to  the proof of part  ($ii$) of Proposition \ref{rcc-se}).

Then we {\bf reformulate} the distributed computation of \eqref{AXBeF} with RRR structure as the following optimization problem
\begin{subequations}\label{opt-rrr}
\begin{align}
\min_{X,Y_{\rm E},Z}&\quad\sum_{i=1}^n\|A_{vi}Y_i -F_{vi}\|_{\rm F}^2, \label{opt-rrr1}\\
\text{s. t. }&\quad \frac{1}{n}Y_i-X_{li}B_{vi} +\sum_{j=1}^n a_{i,j}(Z_i-Z_j) = 0_{r \times q}, \nonumber\\
&\quad Y_i=Y_j,\quad i,j\in\{1,\ldots,n\},\label{opt-rrr2}
\end{align}
\end{subequations}
where $X = [X_{l1},\ldots,X_{ln}]\in\mathbb R^{r\times p}$, $Y_{\rm E}=[Y_1^{\rm T},\cdots, Y_n^{\rm T}]^{\rm T}\in\mathbb R^{nr\times q}$,  and $Z=[Z_1^{\rm T},\cdots, Z_n^{\rm T}]^{\rm T}\in\mathbb R^{nr\times q}$.

\begin{remark}
In \eqref{opt-rrr2},   $Y_i=Y_j$ is a consensus constraint and $\frac{1}{n}Y_i-X_{li}B_{vi} +\sum_{j=1}^n a_{i,j}(Z_i-Z_j) = 0_{r \times q}$ is a coupled constraint, which may be viewed as a (generalized) resource allocation constraint  \cite{YHL:Automatica:2016}.
\hfill$\Diamond$
\end{remark}

It is not hard to obtain the following result.

\begin{proposition}\label{rrr-se}
Suppose that the undirected graph $\mathcal G$ is connected.   $X^*\in\mathbb R^{r\times p}$ is a least squares solution to matrix equation \eqref{AXBeF} if and only if there exist $Y_{\rm E}^*\in\mathbb R^{nr\times q}$ and $Z^*\in\mathbb R^{nr\times q}$ such that $({X}^*,Y_{\rm E}^*, Z^*)$ is a solution to problem \eqref{opt-rrr}.
\end{proposition}

The proof  is omitted due to the space limitation and similarity to that of Proposition \ref{rcc-se}.

In RRR structure, we define $\Lambda^1 = \begin{bmatrix}(\Lambda^1_1)^{\rm T},
\cdots,(\Lambda^1_n)^{\rm T}
\end{bmatrix}^{\rm T}\in\mathbb R^{nr\times q}$ and $\Lambda^2 = \begin{bmatrix}(\Lambda^2_1)^{\rm T},
\cdots, (\Lambda^2_n)^{\rm T}
\end{bmatrix}^{\rm T}\in\mathbb R^{nr\times q}$ as estimates of Lagrangian multipliers, where $\Lambda^1_i\in\mathbb R^{r\times q}$ and $\Lambda^2_i\in\mathbb R^{r\times q}$ for $i\in\{1,\ldots,n\}$.
Then we propose a distributed {\bf algorithm} of agent $i$ as follows:
\begin{subequations}\label{algorithm1}
\begin{align}
\dot X_{li}(t) &= \Lambda^1_i(t) B_{vi}^{\rm T},\quad X_{li}(0)=X_{li0}\in\mathbb R^{r\times p_i},\\
\dot Y_i(t) &= -A_{vi}^{\rm T}(A_{vi}Y_i(t)-F_{vi}) -\sum_{j=1}^n a_{i,j}(Y_i(t)-Y_j(t))\nonumber\\
&\quad -\frac{1}{n}\Lambda^1_i(t)-\sum_{j=1}^n a_{i,j}(\Lambda^2_i(t)-\Lambda^2_j(t)),\nonumber\\
&\quad  Y_i(0)=Y_{i0}\in\mathbb R^{r\times q},\\
\dot Z_i(t) &= -\sum_{j=1}^n a_{i,j}(\Lambda^1_i(t)-\Lambda^1_j(t)),\quad Z_i(0)=Z_{i0}\in\mathbb R^{r\times q},\\
\dot \Lambda^1_i(t) &= \frac{1}{n}(Y_i(t){+\dot Y_i(t)})-(X_{li}(t) +{\dot X_{li}(t)})B_{vi} \\
&\quad+\sum_{j=1}^n a_{i,j}(Z_i(t)-Z_j(t))-\sum_{j=1}^n a_{i,j}(\Lambda^1_i(t)-\Lambda^1_j(t)),\nonumber\\
&\quad \Lambda^1_i(0)=\Lambda^1_{i0}\in\mathbb R^{r\times q},\\
\dot \Lambda^2_i(t) &= \sum_{j=1}^n a_{i,j}(Y_i(t)-Y_j(t))+{\sum_{j=1}^n a_{i,j}(\dot Y_i(t)-\dot Y_j(t))},\nonumber\\
&\quad \Lambda^2_i(0)=\Lambda^2_{i0}\in\mathbb R^{r\times q},
\end{align}
\end{subequations}
where $i\in\{1,\ldots,n\}$, $t\geq 0$,
$X_{li}(t)$, $Y_i(t)$,  and $Z_i(t)$ are the estimates of solutions to problem \eqref{opt-rrr} by agent $i$ at time $t$, and $a_{i,j}$ is the $(i,j)$th element of the adjacency matrix of graph $\mathcal G$.

\begin{remark}\label{derivativeR}
The derivative feedbacks $\dot X_{li}$ and $\dot Y_i$ are used in algorithm \eqref{algorithm1}; otherwise the trajectories of the algorithm may oscillate following a periodic routine. In fact,  derivative feedbacks play a role as a damping term to deal with the general convexity of objective functions \cite{Antipin:2003}.
\hfill $\Diamond$
\end{remark}

The following result shows the correctness  of algorithm \eqref{algorithm1} for problem \eqref{opt-rrr}.

\begin{lemma}\label{rrr-eq}
Suppose that the undirected graph $\mathcal G$ is connected.  $(X^*,Y_{\rm E}^*,Z^*)\in\mathbb R^{r\times p}\times \mathbb R^{nr\times q} \times \mathbb R^{nr\times q}$ is a solution to problem \eqref{opt-rrr} if and only if   there exist matrices $\Lambda^{1*}\in \mathbb R^{nr\times q}$ and $\Lambda^{2*}\in \mathbb R^{nr\times q}$ such that $(X^*,Y_{\rm E}^*,Z^*, \Lambda^{1*},\Lambda^{2*})$  is an equilibrium of \eqref{algorithm1}.
\end{lemma}

The proof is omitted because it is easy due to the KKT optimality condition \cite{Ruszczynski:2006}.

Then we show the convergence of algorithm \eqref{algorithm1}.
Define a function
\begin{eqnarray}\label{V_RRR}
V(X,Y_{\rm E},Z,\Lambda^{1},\Lambda^{2}) = V_1(Y_{\rm E})+ V_2( {X},Y_{\rm E},Z,\Lambda^1,\Lambda^2),
\end{eqnarray}
where
\begin{eqnarray*}
V_1  &\triangleq & \frac{1}{2}\sum_{i=1}^n\|A_{vi}Y_i -F_{vi}\|_{\rm F}^2 -\frac{1}{2} \sum_{i=1}^n\|A_{vi}Y_i^* -F_{vi}\|_{\rm F}^2\\
&&-\sum_{i=1}^n\langle A_{vi}^{\rm T}(A_{vi}Y_i^* -F_{vi}), Y_i -Y_i^* \rangle_{\rm F}\\
&&+ \frac{1}{2}\sum_{i=1}^n\sum_{j=1}^n a_{i,j}\langle Y_i, Y_i -Y_j \rangle_{\rm F},
\end{eqnarray*}
\begin{eqnarray*}
V_2  & \triangleq  &\frac{1}{2}\sum_{i=1}^n \|X_{li}-X_{li}^*\|_{\rm F}^2 +  \frac{1}{2}\sum_{i=1}^n \|Y_i-Y_i^*\|_{\rm F}^2 \\
&&+ \frac{1}{2} \sum_{i=1}^n \| Z_i-Z^{*}_i\|_{\rm F}^2  +  \frac{1}{2}\sum_{i=1}^n \|\Lambda^1_i-\Lambda^{1*}_i\|_{\rm F}^2\\
&&
+ \frac{1}{2} \sum_{i=1}^n \|\Lambda^2_i-\Lambda^{2*}_i\|_{\rm F}^2,
\end{eqnarray*}
and $( X^*,Y_{\rm E}^*,Z^*,\Lambda^{1*},\Lambda^{2*})$  is an equilibrium of \eqref{algorithm1}. The following lemma will be needed in the theoretical proof of our algorithm.

\begin{lemma}\label{V_rrr_p}
If the undirected graph $\mathcal G$ is connected, $V_1(Y_{\rm E})$ defined in \eqref{V_RRR} is nonnegative for all $Y_{\rm E}\in \mathbb R^{nr\times q}$.
\end{lemma}

\begin{IEEEproof}
Consider function $V_1(Y_{\rm E})$ defined in \eqref{V_RRR}.  Because $\mathcal G$ is undirected,
\begin{eqnarray*}
\frac{1}{2}\sum_{i=1}^n\sum_{j=1}^n a_{i,j}\langle Y_i, Y_i -Y_j \rangle_{\rm F}
=\frac{1}{4}\sum_{i=1}^n\sum_{j=1}^n a_{i,j}\| Y_i -Y_j \|^2_{\rm F}\geq 0.
\end{eqnarray*}

Define $ f(Y_{\rm E})=\frac{1}{2}\sum_{i=1}^n\|A_{vi}Y_i -F_{vi}\|_{\rm F}^2$. $f(Y_{\rm E})$ is clearly convex with respect to  matrix $Y_{\rm E}\in\mathbb R^{nr\times q}$.   Then
\begin{eqnarray*}
f(Y_{\rm E})-f(Y_{\rm E}^*) & \geq & \sum_{i=1}^n \langle Y_{i}-Y_{i}^*, \nabla_{Y_{i}^*}f(Y_{\rm E}^*)\rangle_{\rm F}\\
& = & \sum_{i=1}^n \langle Y_{i}-Y_{i}^*, A_{vi}^{\rm T}(A_{vi}Y_i^* -F_{vi})\rangle_{\rm F}.
\end{eqnarray*}
Hence, $V_1(Y_{\rm E})\geq f(Y_{\rm E})-f(Y_{\rm E}^*) -\sum_{i=1}^n \langle Y_{i}-Y_{i}^*, A_{vi}^{\rm T}(A_{vi}Y_i^* -F_{vi})\rangle_{\rm F}\geq 0$ for all $Y_{\rm E}\in \mathbb R^{nr\times q}$.
\end{IEEEproof}

Next, we show the {\bf exponential  convergence} of algorithm \eqref{algorithm1}.
\begin{proposition}\label{Thm-RRR}
If the undirected graph $\mathcal G$ is connected, then
\begin{enumerate}
  \item every equilibrium of algorithm \eqref{algorithm1} is Lyapunov stable and its trajectory is bounded for any initial condition;
  \item {the trajectory of algorithm \eqref{algorithm1} is exponentially convergent and}  $X(t)$ converges to a least squares solution of \eqref{AXBeF} exponentially.
\end{enumerate}
\end{proposition}

The proof can be found in Appendix \ref{app-Thm-RRR}.

\begin{remark}
Compared with related results in linear algebraic equations or others \cite{SA:TAC:2015,MLM:TAC:2015, AMMH:Arxiv:2015,CZH:CCC:2017,LMM:TAC:2017,LMNM:Automatica:2017}, the boundedness assumption for least squares solutions or the existence of exact solutions is not required.\hfill $\Diamond$
\end{remark}

\subsection{Column-Column-Row Structure}

To handle CCR structure, we take a substitutional variable
$$
Y = \begin{bmatrix} Y_{v1}\\
\vdots\\
Y_{vn}
\end{bmatrix}\in\mathbb R^{r\times q},\quad Y_{vi}\in\mathbb R^{r_i\times q},\quad \forall i\in\{1,\ldots,n\}.
$$
It is clear that \eqref{AXBeF}  and \eqref{CCR} is equivalent to $AY= F,\; Y = XB.$  Let $X_i\in\mathbb R^{r\times p}$ be the estimate of $X$ by agent $i$. Define  matrices $[Y_{vi}]^{\{r_j\}_{j=1}^n}_{{\rm R}}$, $[F_{vi}]^{\{m_j\}_{j=1}^n}_{{\rm R}}$, and $[B_{li}]^{\{q_j\}_{j=1}^n}_{{\rm C}}$  as in \eqref{r-ci} and \eqref{r-cii} and  take  $[Y_{vi}]_{\rm R}$, $[F_{vi}]_{\rm R}$, and  $[B_{li}]_{\rm C}$ to represent $[Y_{vi}]^{\{r_j\}_{j=1}^n}_{{\rm R}}$, $[F_{vi}]^{\{m_j\}_{j=1}^n}_{{\rm R}}$, and $[B_{li}]^{\{q_j\}_{j=1}^n}_{{\rm C}}$  for the ease of notation.
Clearly, $\sum_{i=1}^n [Y_{vi}]_{\rm R}=Y$, $\sum_{i=1}^n [F_{vi}]_{\rm R}=F$, and $\sum_{i=1}^n [B_i]_{\rm C}=B$.    Here we construct a transformation, called {\bf $X$-consensus substitutional decomposition} with requiring the consensus of $X_i$,  and then \eqref{AXBeF}  and \eqref{CCR} is equivalent to
\begin{align}
\sum_{i=1}^n A_{li}Y_{vi} =& \sum_{i=1}^n [F_{vi}]_{\rm R}, \label{CCR-DP1}\\
\sum_{i=1}^n [Y_{vi}]_{\rm R} =& \sum_{i=1}^n X_i[B_{li}]_{\rm C},\, X_i=X_j,\, i,j\in\{1,\ldots,n\}. \label{CCR-DP2}
\end{align}

To decompose \eqref{CCR-DP1} and \eqref{CCR-DP2}, we add new variables $U_i\in\mathbb R^{m\times q}$,  $W_i\in\mathbb R^{m\times q}$, and $Z_i\in\mathbb R^{r\times q}$ such that
\begin{align}
& A_{li}Y_{vi} - [F_{vi}]_{\rm R}-U_i=0_{m\times q}, \, U_i = \sum_{j=1}^n a_{i,j}(W_i-W_j),\label{CCR-DD1}\\
&  [Y_{vi}]_{\rm R} - X_i[B_{li}]_{\rm C}-\sum_{j=1}^n a_{i,j}(Z_i-Z_j)= 0_{r\times q},\, X_i=X_j,\label{CCR-DD2}
\end{align}
where $i,j\in\{1,\ldots,n\},$ $a_{i,j}$ is the $(i,j)$th element of the adjacency matrix of graph $\mathcal G$. If \eqref{CCR-DD1} and \eqref{CCR-DD2} hold, then one can easily obtain \eqref{CCR-DP1} and \eqref{CCR-DP2}. Conversely, if \eqref{CCR-DP1} and \eqref{CCR-DP2} hold, it follows from a similar proof of Proposition \ref{rcc-se} that there exist $U_i$, $W_i$, and $Z_i$ such that \eqref{CCR-DD1} and \eqref{CCR-DD2} hold.

Let $X_{\rm E} = [X_1^{\rm T},\ldots,X_n^{\rm T}]^{\rm T}\in\mathbb R^{r\times p}$, $Y = [Y_{v1}^{\rm T},
\cdots,
Y_{vn}^{\rm T}
]^{\rm T}\in\mathbb R^{r\times q}$, $U=[ U_1^{\rm T},
\cdots,
U_n^{\rm T}
]^{\rm T}\in\mathbb R^{nm\times q}$, $W=[ W_1^{\rm T},
\cdots, W_n^{\rm T} ]^{\rm T}\in\mathbb R^{nm\times q}$, and $Z=[Z_1^{\rm T},
\cdots, Z_n^{\rm T} ]^{\rm T}\in\mathbb R^{nr\times q}$.
Then we {\bf reformulate} the distributed computation of \eqref{AXBeF} with CCR structure as the following optimization problem
\begin{subequations}\label{opt-ccr}
\begin{align}
\min_{X_{\rm E},Y,U,W,Z}&\quad\sum_{i=1}^n\|A_{li}Y_{vi} - [F_{vi}]_{\rm R}-U_i\|_{\rm F}^2, \label{opt-ccr1}\\
\text{s. t.}&\quad X_i=X_j,\, U_i = \sum_{j=1}^n a_{i,j}(W_i-W_j),\\
 &\, [Y_{vi}]_{\rm R} - X_i[B_{li}]_{\rm C}-\sum_{j=1}^n a_{i,j}(Z_i-Z_j)= 0_{r\times q},\label{opt-ccr2}\\
 &\, i,j\in\{1,\ldots,n\}.\nonumber
\end{align}
\end{subequations}

Similar to RRR structure,  \eqref{CCR-DP1} and \eqref{CCR-DP2}  are a combination of  a consensus constraint and coupled equality constraints.
Then we have the following result.

\begin{proposition}\label{ccr-se}
Suppose that the undirected graph $\mathcal G$ is connected.  $X^*\in\mathbb R^{r\times p}$ is a least squares solution to matrix equation \eqref{AXBeF} if and only if  there exist $X_{\rm E}^* = 1_n\otimes X^*$, $Y^*\in\mathbb R^{r\times q}$, $Z^*\in\mathbb R^{nr\times q}$, $U^*\in\mathbb R^{nm\times q}$, and $W^*\in\mathbb R^{nm\times q}$ such that $(X_{\rm E}^*,Y^*,Z^*,U^*,W^*)$ is a solution to problem \eqref{opt-ccr}.
\end{proposition}

The proof  is omitted due to the space limitation and similarity to that of Proposition \ref{rcc-se}.

In  CCR structure, we take $\Lambda^1 = \begin{bmatrix}\Lambda^1_1\\
\vdots\\
\Lambda^1_n
\end{bmatrix}\in\mathbb R^{nr\times p}$, $\Lambda^2 = \begin{bmatrix}\Lambda^2_1\\
\vdots\\
\Lambda^2_n
\end{bmatrix}\in\mathbb R^{nm\times q}$, and $\Lambda^3 = \begin{bmatrix}\Lambda^3_1\\
\vdots\\
\Lambda^3_n
\end{bmatrix}\in\mathbb R^{nr\times q}$ as the Lagrangian multipliers, where $\Lambda^1_i\in\mathbb R^{r\times p}$, $\Lambda^2_i\in\mathbb R^{m\times q}$, and $\Lambda^3_i\in\mathbb R^{r\times q}$.
The distributed {\bf algorithm} of agent $i$ is
\begin{subequations}\label{algorithm3}
\begin{align}
\dot X_i(t) &= \Lambda^3_i(t)[B_{li}]_{\rm C}^{\rm T}-\sum_{j=1}^n a_{i,j}(\Lambda^1_i(t)-\Lambda^1_j(t))\nonumber\\
&\quad -\sum_{j=1}^n a_{i,j}(X_i(t)-X_j(t)),\, X_i(0)=X_{i0}\in\mathbb R^{r\times p},\\
\dot Y_{vi}(t) &= -A_{li}^{\rm T}(A_{li}Y_{vi}(t) - [F_{vi}]_{\rm R}-U_i(t))-[I_{r_i}]_{\rm C}\Lambda^3_i(t),\nonumber\\
&\quad Y_{vi}(0)=Y_{vi0}\in\mathbb R^{r_i\times q},\\
\dot U_i(t) &= A_{li}Y_{vi}(t) - [F_{vi}]_{\rm R}-U_i(t)-\Lambda^2_i(t),\nonumber\\
&\quad U_i(0)=U_{i0}\in\mathbb R^{m\times q},\\
\dot W_i(t) &= \sum_{j=1}^n a_{i,j}(\Lambda^2_i(t)-\Lambda^2_j(t)),\, W_i(0)=W_{i0}\in\mathbb R^{m\times q},\\
\dot Z_i(t) &= \sum_{j=1}^n a_{i,j}(\Lambda^3_i(t)-\Lambda^3_j(t)),\, Z_i(0)=Z_{i0}\in\mathbb R^{r\times q},\\
\dot \Lambda^1_i(t) &=  \sum_{j=1}^n a_{i,j}(X_i(t)-X_j(t)),\quad \Lambda^1_i(0)=\Lambda^1_{i0}\in\mathbb R^{r\times p},\\
\dot \Lambda^2_i(t) &=  U_i(t)+{ \dot U_i(t)}- \sum_{j=1}^n a_{i,j}(W_i(t)-W_j(t))\nonumber\\
&\,-\sum_{j=1}^n a_{i,j}(\Lambda^2_i(t)-\Lambda^2_j(t)),\, \Lambda^2_i(0)=\Lambda^2_{i0}\in\mathbb R^{m\times q},\\
\dot \Lambda^3_i(t) &= [Y_{vi}]_{\rm R}(t)+ {[\dot Y_{vi}]_{\rm R}(t)} - X_{i}(t)[B_{li}]_{\rm C}\nonumber\\
&\quad-\sum_{j=1}^n a_{i,j}(Z_i(t)-Z_j(t))-\sum_{j=1}^n a_{i,j}(\Lambda^3_i(t)-\Lambda^3_j(t)),\nonumber\\
&\quad  \Lambda^3_i(0)=\Lambda^3_{i0}\in\mathbb R^{r\times q},
\end{align}
\end{subequations}
where $i\in\{1,\ldots,n\}$, $t\geq 0$,
$X_i(t)$, $Y_{vi}(t)$, $U_i(t)$,  $W_i(t)$, and $Z_i(t)$ are the estimates of solutions to problem \eqref{opt-ccr} by agent $i$ at time $t$, $a_{i,j}$ is the $(i,j)$th element of the adjacency matrix of graph $\mathcal G$, and $[B_{li}]_{\rm C}$, $[F_{vi}]_{\rm R}$, and $[I_{r_i}]_{\rm C}$ are shorthand notations for $[B_{li}]^{\{q_j\}_{j=1}^n}_{{\rm C}}$, $[F_{vi}]^{\{m_j\}_{j=1}^n}_{{\rm R}}$, and $[I_{r_i}]^{\{r_j\}_{j=1}^n}_{{\rm C}}$ as defined in \eqref{r-ci} and \eqref{r-cii}.

Similar to algorithm \eqref{algorithm1}, algorithm \eqref{algorithm3} is the saddle-point dynamics of the modified Lagrangian function with derivative feedbacks, which are a ``damping" term (see Remark \ref{derivativeR}).

The following lemma reveals the connection of solutions to problem \eqref{opt-ccr} and equilibria of algorithm \eqref{algorithm3},
whose proof is quite obvious because of the KKT optimality condition   \cite{Ruszczynski:2006}.

\begin{lemma}\label{ccr-eq}
Suppose that the undirected graph $\mathcal G$ is connected.   $(X_{\rm E}^*,Y^*,Z^*,U^*,W^*)\in\mathbb R^{nr\times p}\times \mathbb R^{r\times q}\times \mathbb R^{nr\times q} \times \mathbb R^{nm\times q}\times \mathbb R^{nm\times q}$ is a solution to problem \eqref{opt-ccr} if and only if  there exist matrices $\Lambda^{1*}\in \mathbb R^{nr\times p}$, $\Lambda^{2*}\in \mathbb R^{nm\times q}$, and $\Lambda^{3*}\in \mathbb R^{nr\times q}$ such that $(X_{\rm E}^*,Y^*,Z^*,U^*,W^*,\Lambda^{1*},\Lambda^{2*},\Lambda^{3*})$  is an equilibrium of \eqref{algorithm3}.
\end{lemma}

Define the function
\begin{align}\label{V_CCR}
V(X_{\rm E},Y,Z,U,W,\Lambda^{1},\Lambda^{2},\Lambda^{3}) =& V_1(X_{\rm E},Y,U)\nonumber\\
&\hspace{-2 cm}+ V_2(X_{\rm E},Y,Z,U,W,\Lambda^{1},\Lambda^{2},\Lambda^{3})
\end{align}
with
\begin{align*}
V_1 \triangleq & \frac{1}{2}\sum_{i=1}^n\|A_{li}Y_{vi} - [F_{vi}]_{\rm R}-U_i\|_{\rm F}^2 \nonumber\\
 &+\sum_{i=1}^n\sum_{j=1}^n\langle \Lambda^{1*}_i,  a_{i,j}(X_i-X_j)\rangle_{\rm F}\nonumber\\
&+\sum_{i=1}^n\langle \Lambda^{2*}_i, U_i\rangle_{\rm F} +\sum_{i=1}^n\langle \Lambda^{3*}_i, [Y_{vi}]_{\rm R} - X_i[B_{li}]_{\rm C}\rangle_{\rm F}\nonumber\\
& -\frac{1}{2}\sum_{i=1}^n\|A_{li}Y_{vi}^* - [F_{vi}]_{\rm R}-U_i^*\|_{\rm F}^2,
\end{align*}
\begin{align*}
V_2 \triangleq & \frac{1}{2}\sum_{i=1}^n \|X_{i}-X_{i}^*\|_{\rm F}^2+  \frac{1}{2}\sum_{i=1}^n \|Y_{vi}-Y_{vi}^*\|_{\rm F}^2 \nonumber\\
 &+  \frac{1}{2}\sum_{i=1}^n \|Z_{i}-Z_{i}^*\|_{\rm F}^2 +  \frac{1}{2}\sum_{i=1}^n \|U_{i}-U_{i}^*\|_{\rm F}^2 \\
&+   \frac{1}{2}\sum_{i=1}^n \|W_{i}-W_{i}^*\|_{\rm F}^2  +  \frac{1}{2}\sum_{i=1}^n \|\Lambda^1_i-\Lambda^{1*}_i\|_{\rm F}^2\nonumber\\
 &  + \frac{1}{2} \sum_{i=1}^n \|\Lambda^2_i-\Lambda^{2*}_i\|_{\rm F}^2+ \frac{1}{2} \sum_{i=1}^n \|\Lambda^3_i-\Lambda^{3*}_i\|_{\rm F}^2,
\end{align*}
where $(X_{\rm E}^*,Y^*,Z^*,U^*,W^*,\Lambda^{1*},\Lambda^{2*},\Lambda^{3*})$  is an equilibrium of \eqref{algorithm3}.

\begin{lemma}\label{V_1_p}
Suppose that the undirected graph $\mathcal G$ is connected. The function $V_1(X_{\rm E},Y,U)$ defined in \eqref{V_CCR} is nonnegative for all $(X_{\rm E},Y,U)\in\mathbb R^{nr\times p}\times \mathbb R^{r\times q}\times \mathbb R^{nm\times q}$.
\end{lemma}

{The proof is similar to that of Lemma \ref{V_rrr_p} and omitted.}

The following result shows the {\bf exponential convergence} of algorithm \eqref{algorithm3}.
\begin{proposition}\label{THM3}
If the undirected graph $\mathcal G$ is connected, then
\begin{enumerate}
  \item every equilibrium of algorithm \eqref{algorithm3} is Lyapunov stable and its trajectory is bounded for any initial condition;
  \item {the trajectory of algorithm \eqref{algorithm3} is exponentially convergent and} $X_i(t)$ converges to a least squares solution of \eqref{AXBeF} exponentially for all $i\in\{1,\ldots,n\}$.
\end{enumerate}
\end{proposition}

The proof can be found in Appendix \ref{app-THM3}.

\subsection{Column-Row-Row Structure}

In   CRR structure, which is the most complicated structure among the 4 standard structures, the above decomposition methods do not work.
Define a substitutional variable
$$
Y = \begin{bmatrix} Y_{v1}\\
\vdots\\
Y_{vn}
\end{bmatrix}\in\mathbb R^{r\times q}, \quad Y_{vi}\in\mathbb R^{r_i\times q},\quad \forall i\in\{1,\ldots,n\}.
$$
Clearly, \eqref{AXBeF} with
\eqref{CRR} is equivalent to
$AY= F$ and $ Y = XB.$
Moreover, we further define the augmented matrices $[Y_{vi}]^{\{r_j\}_{j=1}^n}_{{\rm R}}$ and $[F_{vi}]^{\{m_j\}_{j=1}^n}_{{\rm R}}$ as  in \eqref{r-ci} and take $[Y_{vi}]_{\rm R} $ and $[F_{vi}]_{\rm R}$ to denote $[Y_{vi}]^{\{r_j\}_{j=1}^n}_{{\rm R}}$ and $[F_{vi}]^{\{m_j\}_{j=1}^n}_{{\rm R}}$ for convenience.
Then we have
\begin{eqnarray}
\sum_{i=1}^n A_{li}Y_{vi} &=& \sum_{i=1}^n [F_{vi}]_{\rm R} ,\label{CRR-DP1} \\
\sum_{i=1}^n [Y_{vi}]_{\rm R}  &=& \sum_{i=1}^n X_{li}B_{vi}.\label{CRR-DP2}
\end{eqnarray}

To decompose \eqref{CRR-DP1} and \eqref{CRR-DP2}, we take new variables $U_i\in\mathbb R^{m\times q}$,  $W_i\in\mathbb R^{m\times q}$ and $Z_i\in\mathbb R^{r\times q}$ such that
\begin{align}
& A_{li}Y_{vi} - [F_{vi}]_{\rm R}-U_i=0_{m\times q}, \, U_i = \sum_{j=1}^n a_{i,j}(W_i-W_j),\label{CRR-DD1}\\
&  [Y_{vi}]_{\rm R} - X_{li}B_{vi}-\sum_{j=1}^n a_{i,j}(Z_i-Z_j)= 0_{r\times q}. \label{CRR-DD2}
\end{align}

Let $X = [X_{l1},\ldots,X_{ln}]\in\mathbb R^{r\times p}$,
$Y = [ Y_{v1}^{\rm T},\cdots, Y_{vn}^{\rm T}]^{\rm T}\in\mathbb R^{r\times q}$,
$U = [U_1^{\rm T},\ldots,U_n^{\rm T}]^{\rm T}\in\mathbb R^{nm\times q}$,   $W = [W_1^{\rm T},\ldots,W_n^{\rm T}]^{\rm T}\in\mathbb R^{nm\times q}$,  and $Z = [Z_1^{\rm T},\ldots,Z_n^{\rm T}]^{\rm T}\in\mathbb R^{nr\times q}$. We {\bf reformulate} the distributed computation of \eqref{AXBeF} with CRR structure as  the following optimization problem
\begin{subequations}\label{opt-crr}
\begin{align}
\min_{X,Y, U, W,Z}&\quad \sum_{i=1}^n\|A_{li}Y_{vi} - [F_{vi}]_{\rm R}-U_i\|_{\rm F}^2, \label{opt-crr1}\\
\text{s. t.}&\quad  \,[Y_{vi}]_{\rm R} - X_{li}B_{vi}- \sum_{j=1}^n a_{i,j}(Z_i-Z_j)= 0_{r\times q},\label{opt-crr2}\\
&\quad U_i = \sum_{j=1}^n a_{i,j}(W_i-W_j),\quad  i\in\{1,\ldots,n\}.
\end{align}
\end{subequations}

The transformation given here is simply called {\bf consensus-free substitutional decomposition} because we do not need the consensus of $X_i$ or $Y_i$ for $i=1,\ldots,n$.  Then we have the following result.

\begin{proposition}\label{crr-se}
Suppose that the undirected graph $\mathcal G$ is connected.  $X^*$ is a least squares solution to equation \eqref{AXBeF} if and only if there exist $Y^*$, $Z^*$, $U^*$, and $W^*$ such that $(X^*,Y^*,Z^*,U^*,W^*)$ is a solution to problem \eqref{opt-crr}.
\end{proposition}

The proof  is omitted due to the space limitation and similarity to that of Proposition \ref{rcc-se}.

In  this structure, we propose a distributed {\bf algorithm} of agent $i$ as follows:
\begin{subequations}\label{algorithm4}
\begin{align}
\dot X_{li}(t)&= \Lambda^2_i (t) B_{vi}^{\rm T},\quad X_{li}(0)= X_{li0}\in\mathbb R^{r\times p_i},\\
\dot Y_{vi} (t)&= -A_{li}^{\rm T}(A_{li}Y_{vi}(t) - [F_{vi}]_{\rm R}-U_i(t))-[I_{r_i}]_{\rm C}\Lambda^2_i(t),\nonumber\\
&\quad Y_{vi} (0)=Y_{vi0} \in\mathbb R^{r_i\times q}\\
\dot U_i (t)&= A_{li}Y_{vi}(t) - [F_{vi}]_{\rm R}-U_i(t)-\Lambda^1_i(t),\nonumber\\
&\quad U_i (0)=U_{i0}\in\mathbb R^{m\times q},\\
\dot W_i (t) &= \sum_{j=1}^n a_{i,j}(\Lambda^1_i (t) -\Lambda^1_j(t)),\, W_i (0)=W_{i0}\in\mathbb R^{m\times q},\\
\dot Z_i (t)&= \sum_{j=1}^n a_{i,j}(\Lambda^2_i(t)-\Lambda^2_j(t)),\quad Z_i (0)=Z_{i0}\in\mathbb R^{r\times q},\\
\dot \Lambda^1_i (t)&=  U_i(t)+{\dot U_i(t)}- \sum_{j=1}^n a_{i,j}(W_i(t)-W_j(t))\nonumber\\
&\quad-\sum_{j=1}^n a_{i,j}(\Lambda^1_i(t)-\Lambda^1_j(t)),\, \Lambda^1_i (0)=\Lambda^1_{i0}\in\mathbb R^{m\times q}, \\
\dot \Lambda^2_i(t) &=  [Y_{vi}]_{\rm R}(t) +{ [\dot Y_{vi}]_{\rm R}}(t)- X_{li}(t)B_{vi}\nonumber\\
&\quad-\sum_{j=1}^n a_{i,j}(Z_i(t)-Z_j(t)) -\sum_{j=1}^n a_{i,j}(\Lambda^2_i(t)-\Lambda^2_j(t))\nonumber\\
&\quad- \dot X_{li}(t)B_{vi},
\quad \Lambda^2_i (0)=\Lambda^2_{i0}\in\mathbb R^{r\times q},
\end{align}
\end{subequations}
where $i\in\{1,\ldots,n\}$, $t\geq 0$,
$X_{li}(t)$, $Y_{vi}(t)$, $U_i(t)$,  $W_i(t)$, and $Z_i(t)$ are the estimates of solutions to problem \eqref{opt-crr} by agent $i$ at time $t$,
$a_{i,j}$ is the $(i,j)$th element of the adjacency matrix of graph $\mathcal G$, and
$[Y_{vi}]_{\rm R}=[Y_{vi}]^{\{m_j\}_{j=1}^n}_{{\rm R}}$ and $[I_{r_i}]_{\rm C}=[I_{r_i}]^{\{r_j\}_{j=1}^n}_{{\rm C}}$ are as defined in \eqref{r-ci} and \eqref{r-cii}.
Similar to algorithms \eqref{algorithm1} and \eqref{algorithm3}, the design of algorithm \eqref{algorithm4} also combines the saddle-point dynamics of the modified Lagrangian function and derivative feedback technique.

Let $\Lambda^1 = \begin{bmatrix}\Lambda^1_1\\
\vdots\\
\Lambda^1_n
\end{bmatrix}\in\mathbb R^{nm\times q}$ and $\Lambda^2 = \begin{bmatrix}\Lambda^2_1\\
\vdots\\
\Lambda^2_n
\end{bmatrix}\in\mathbb R^{nr\times q}$, where $\Lambda^1_i\in\mathbb R^{m\times q}$ and $\Lambda^2_i\in\mathbb R^{r\times q}$ for $i\in\{1,\ldots,n\}$. We have the following result, whose proof is omitted because it is straightforward due to the KKT optimality condition (Theorem 3.25 of \cite{Ruszczynski:2006}).

\begin{lemma}\label{crr-eq}
Suppose that the undirected graph $\mathcal G$ is connected.  $(X^*,Y^*,Z^*,U^*,W^*)\in\mathbb R^{r\times p} \times \mathbb R^{r\times q} \times \mathbb R^{nr\times q} \times \mathbb R^{nm\times q} \times \mathbb R^{nm\times q} $ is a solution to problem \eqref{opt-crr} if and only if   there exist  $\Lambda^{1*}\in\mathbb R^{nm\times q}$ and $\Lambda^{2*}\in\mathbb R^{nr\times q}$ such that $(X^*,Y^*,Z^*,U^*,W^*,\Lambda^{1*},\Lambda^{2*})$  is an equilibrium of \eqref{algorithm4}.
\end{lemma}

For further analysis, let $(X^*,Y^*,Z^*,U^*,W^*,\Lambda^{1*},\Lambda^{2*})$  be an equilibrium of \eqref{algorithm4},
and take
\begin{eqnarray}\label{V_CRR}
V(X,Y,Z,U,W,\Lambda^{1},\Lambda^{2}) &=& V_1(Y,U)\nonumber\\
&&\hspace{-1 cm}+ V_2(X,Y,Z,U,W,\Lambda^{1},\Lambda^{2}),
\end{eqnarray}
where
\begin{eqnarray*}
V_1 &=&  \frac{1}{2}\sum_{i=1}^n\|A_{li}Y_{vi} - [F_{vi}]_{\rm R}-U_i\|_{\rm F}^2+\sum_{i=1}^n\langle \Lambda^{1*}_i, U_i\rangle_{\rm F}\nonumber\\
&& +\sum_{i=1}^n\langle \Lambda^{2*}_i, [Y_{vi}]_{\rm R} \rangle_{\rm F} - \frac{1}{2}\sum_{i=1}^n\|A_{li}Y_{vi}^* - [F_{vi}]_{\rm R}-U_i^*\|_{\rm F}^2
\end{eqnarray*}
and
\begin{eqnarray*}
V_2 &\triangleq & \frac{1}{2}\sum_{i=1}^n \|X_{li}-X_{li}^*\|_{\rm F}^2 +  \frac{1}{2}\sum_{i=1}^n \|Y_{vi}-Y_{vi}^*\|_{\rm F}^2\\
&& +  \frac{1}{2}\sum_{i=1}^n \|Z_{i}-Z_{i}^*\|_{\rm F}^2 +  \frac{1}{2}\sum_{i=1}^n \|U_{i}-U_{i}^*\|_{\rm F}^2    \\
&&+   \frac{1}{2}\sum_{i=1}^n \|W_{i}-W_{i}^*\|_{\rm F}^2+  \frac{1}{2}\sum_{i=1}^n \|\Lambda^1_i-\Lambda^{1*}_i\|_{\rm F}^2 \\
&&+ \frac{1}{2} \sum_{i=1}^n \|\Lambda^2_i-\Lambda^{2*}_i\|_{\rm F}^2.
\end{eqnarray*}
Then we get the following result.

\begin{lemma}\label{V_2_p}
The function $V_1(  Y,U)$ defined in \eqref{V_CRR} is nonnegative for all $(Y,U)\in\mathbb R^{r\times q}\times \mathbb R^{nm\times q}$.
\end{lemma}

The proof is similar to that of Lemma \ref{V_1_p} and is omitted.
Now it is time to show the main result of this subsection.

Next, we show the {\bf exponential convergence} of algorithm \eqref{algorithm4}.
\begin{proposition}\label{THM-CRR}
If the undirected graph $\mathcal G$ is connected, then
\begin{enumerate}
  \item every equilibrium of algorithm \eqref{algorithm4} is Lyapunov stable and its trajectory is bounded for any initial condition;
  \item {the trajectory of algorithm \eqref{algorithm4} is exponentially convergent and} $X(t)$ converges to a least squares solution of the matrix equation  \eqref{AXBeF} exponentially.
\end{enumerate}
\end{proposition}

The proof can be found in Appendix \ref{app-THM-CRR}.

\subsection{Discussions}
The conclusion of Theorem \ref{themain} is obtained immediately from  the results given in Propositions \ref{rcc-se}-\ref{THM-CRR}.   In fact, we develop new methods for the distributed computation of a least squares solution to matrix equation \eqref{AXBeF}, which is much more complicated than that to the linear algebraic equation. The main results of this section is summarized as follows:
\begin{itemize}
\item  We employ different substitutional decomposition methods to {\bf reformulate} the original computation matrix equations as distributed constrained optimization problems with different constraints in the standard structures.  Note that the decompositions are new compared with those in the distributed computation of the linear algebraic equation of the form $Ax=b$  in \cite{MLM:TAC:2015,SA:TAC:2015,WFM:ACC:2016,SA:ACC:2016,LMM:TAC:2017,LMNM:Automatica:2017}.
\item  We give distributed  {\bf algorithms} to deal with the distributed constrained optimization problems, which are equivalent to matrix equations in different standard structures.    { The proposed algorithms are not a direct application of the existing ideas on distributed subgradient optimization designs. Derivative feedback ideas are used to deal with the general convexity. Additionally, the auxiliary variables are employed as observers to estimate the unavailable matrix information from the structure decomposition. Therefore, the proposed algorithms are different from those given in  \cite{KCM:A:2015} and \cite{SA:TAC:2015}, which did not used derivative feedbacks and auxiliary variables.}
\item  We give the {\bf exponential convergence} analysis  of the algorithms by using advanced (control) techniques such as the Lyapunov stability theory and the derivative feedback  to deal with convexity of objective functions. The proposed algorithms are globally exponentially convergent, which guarantees the exponential convergence  to a least squares solution to the matrix equation for any initial condition.
\end{itemize}

In each standard structure of our problems, we have to employ different ideas to obtain a solution of the reformulated distributed optimization problems because the distributed design for problems with various constraints and only convex objective functions is a nontrivial task.   Moreover, the derivative feedback plays a ``damping" role in the structures with the coupled constraints for the convergence of the proposed algorithms. Specifically,  different consensus variables and derivative feedback variables are used for various structures due to different constraints (see Table \ref{Tab-b}).
The developed approach may provide effective tools for general cases or mixed structures even though there may be no universal way of generalizing this approach.

\begin{table}[H]
\caption{Consensus variable, coupled constraint, and derivative feedback} \label{Tab-b}
\begin{center}
  \begin{tabular}{ l | l | l | l}
    \hline
     &  consensus   & coupled  & derivative \\
          & variable  &  constraint &  feedback\\ \hline
    RCC   & $X_i=X_j$ & none & none\\
       & $Y_i=Y_j$ &  & \\
       \hline

    RRR  & $Y_i=Y_j$ & $\frac{1}{n}\sum_{i=1}^{n} Y_i  = \sum_{i=1}^{n}X_{li}B_{vi}$  & $\dot Y_i$, $\dot X_{li}$\\

  & $\Lambda^1_i=\Lambda^1_j$ &  &  \\ \hline

    CCR    & $X_i=X_j$  & $\sum_{i=1}^n A_{li}Y_{vi} = \sum_{i=1}^n [F_{vi}]_{\rm R}$ & $\dot U_i$, $\dot Y_{vi}$\\
    & $\Lambda^2_i=\Lambda^2_j$ &        $\sum_{i=1}^n [Y_{vi}]_{\rm R} = \sum_{i=1}^n X_i[B_{li}]_{\rm C}$ &\\
       & $\Lambda^3_i=\Lambda^3_j$ &         &\\ \hline

    CRR   &  $\Lambda^1_i=\Lambda^1_j$   & $\sum_{i=1}^n A_{li}Y_{vi} = \sum_{i=1}^n [F_{vi}]_{\rm R}$
             & $\dot U_i$, $\dot Y_{vi}$, $\dot X_{li}$\\
             & $\Lambda^2_i=\Lambda^2_j$   &
            $\sum_{i=1}^n [Y_{vi}]_{\rm R}  = \sum_{i=1}^n X_{li}B_{vi}$  &  \\
    \hline
  \end{tabular}
\end{center}
\end{table}

\begin{remark}
This paper sheds light on  state-of-the-art of  the distributed computation of matrix equations  optimization via a distributed optimization perspective.  For different problem structures, distributed computation algorithms with exponential convergence rates are proposed by combining distributed optimization and control ideas.
\hfill$\Diamond$
\end{remark}

\section{Numerical Simulation}
\label{sec:num}
In this section, we give a numerical example for illustration. Due to the space limitation, we only present a numerical simulation for RRR structure.

Consider a linear matrix equation \eqref{AXBeF} with the structure \eqref{RRR} and $n=4$, where
\begin{eqnarray*}
A_{v1} = [2,\,1],\quad A_{v2} = [4,\,3],\quad A_{v3} = [1,\,3],\quad A_{v4} = [2,\,4],\\
B_{v1} = [1,\,2], \quad B_{v2} = [3,\,2],\quad B_{v3} = [2,\,4],\quad B_{v4} = [2,\,1],
\end{eqnarray*}
and $F$ is given by $$F_{v1} = [0,\,0],\quad F_{v2} = [2,\,1],\quad F_{v3} = [3,\,5],\quad F_{v4} = [1,\,4].$$

\begin{figure}
  \centering
  \includegraphics[width=7 cm]{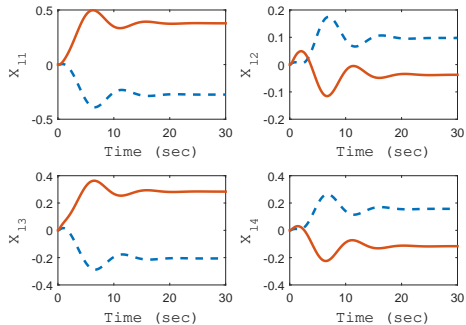}
  \caption{Trajectories of estimates for $X$ versus time}
  \label{X-l}
\end{figure}

\begin{figure}
  \centering
  \includegraphics[width=7 cm]{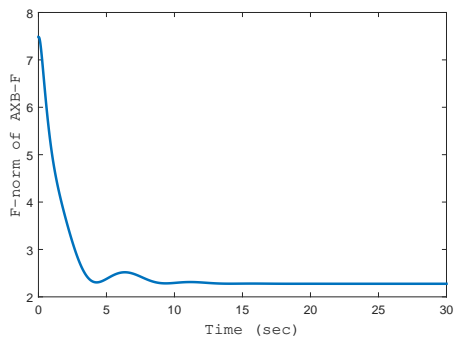}
  \caption{Trajectories of estimates for $\|AXB-F\|_{\rm F}$ versus time}
  \label{FF-l}
\end{figure}

\begin{figure}
  \centering
  \includegraphics[width=7 cm]{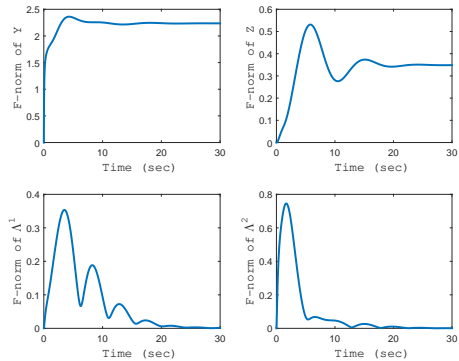}
  \caption{Trajectories of estimates for $Y$, $Z$, $\Lambda^1$, and $\Lambda^2$ versus time }
  \label{auxi-l}
\end{figure}

There is no exact solution for this matrix equation, and therefore, we find a least squares  solution for the problem. Let the adjacency matrix of the graph be
$\begin{bmatrix}0 & 1 & 0 & 1\\
1 & 0 & 1 & 0\\
0 & 1 & 0 & 1\\
1 & 0 & 1 & 0
\end{bmatrix}$.
We solve a least squares solution with algorithm \eqref{algorithm1}
\begin{eqnarray*}
X &=& [X_{l1}, X_{l2}, X_{l3}, X_{l4}]\\
&=&\begin{bmatrix}
-0.2744  &  0.0973  & -0.2058   & 0.1572\\
    0.3780  & -0.0373  &  0.2835  & -0.1163
\end{bmatrix}\in\mathbb R^{2\times 4},
\end{eqnarray*}
where agent $i$ estimates $X_{li}$ for $i\in\{1,\ldots,4\}$. Fig. \ref{X-l} shows that the trajectory of algorithm converges to a least squares solution and Fig. \ref{FF-l} shows the trajectory of $\|AXB-F\|_{\rm F}$, while Fig. \ref{auxi-l} demonstrates the boundedness of algorithm variables.

\section{Conclusions}
\label{conclusion}

In this paper, the distributed computation of least squares solutions to the linear matrix equation $AXB=F$ in standard distributed structures has been studied. Based on substitutional decompositions, the computation problems have been reformulated as equivalent constrained optimization problems in the standard cases. Inspired by saddle-point dynamics and derivative feedbacks, distributed continuous-time algorithms for the reformulated problems have been proposed. Furthermore, the boundedness and exponential convergence of the proposed algorithms have been proved using the  stability and Lyapunov approaches.  Finally, the algorithm performance has  been illustrated via a numerical simulation.

This paper assumes that information of matrices is divided with respect to rows and columns and solves the least squares solutions. It is desirable to further investigate the distributed computation for other well-known linear  matrix equations with mixed structures and solve general solutions such as solutions to  LASSO type problems. In addition,
undirected graphs may be generalized to directed graphs, random graphs, and switching graphs, for instance, and also various effective discrete-time algorithms based on ADMM or other methods may be constructed.

\begin{appendices}\label{App}

\section{Proof of Proposition \ref{rcc-se}}\label{app-rcc-se}

($i$) Suppose that $(X_{\rm E}^*,Y_{\rm E}^*)=(1_n\otimes X^*,1_n\otimes Y^*)$  is a solution to \eqref{opt-rcc}. We show that $X^*$ is a least squares solution to \eqref{AXBeF}.

Because $\mathcal G$ is undirected and connected, \eqref{opt-rcc2} is equivalent to
\begin{eqnarray*}
&& \sum_{j=1}^n a_{i,j}(X_i-X_j)=0_{r\times p},\quad \sum_{j=1}^n a_{i,j}(Y_i-Y_j)=0_{m\times p},\\
&& A_{vi}X_i  = Y_i^{vi},\,i\in\{1,\ldots,n\}.
\end{eqnarray*}

By the KKT optimality condition (Theorem 3.25 of \cite{Ruszczynski:2006}),
$(X_{\rm E}^*,Y_{\rm E}^*)=(1_n\otimes X^*,1_n\otimes Y^*)$  is a solution to problem \eqref{opt-rcc} if and only if
$AX^*=Y^*$ and there are matrices $\Lambda_i^{1*}\in\mathbb R^{r\times p}$, $\Lambda_i^{2*}\in\mathbb R^{m\times p}$, and $\Lambda_i^{3*}\in\mathbb R^{m_i\times p}$ such that
\begin{subequations}\label{opt-rcc-c}
\begin{align}
0_{r\times p} &= -A_{vi}^{\rm T}\Lambda^{3*}_i- \sum_{j=1}^n a_{j,i}(\Lambda^{1*}_i-\Lambda^{1*}_j),\label{opt-rcc-c1}\\
0_{m\times p} &= -(Y^*B_{li} - F_{li})B_{li}^{\rm T}+{[I_{m_i}]_{\rm R} \Lambda^{3*}_i} -\sum_{j=1}^n a_{j,i}(\Lambda^{2*}_i-\Lambda^{2*}_j),\label{opt-rcc-c2}
\end{align}
\end{subequations}
where, for simplicity, $[I_{m_i}]_{\rm R}$ denotes $[I_{m_i}]^{\{m_j\}_{j=1}^n}_{{\rm R}}$ as defined in \eqref{r-ci}.

By \eqref{opt-rcc-c} and $a_{i,j}=a_{j,i}$  because $\mathcal G$ is undirected,  we have
\begin{eqnarray}
0_{r\times p}&=&\sum_{i=1}^{n} [A_{vi}^{\rm T}\Lambda^{3*}_i+ \sum_{j=1}^n a_{j,i}(\Lambda^{1*}_i-\Lambda^{1*}_j)]\nonumber\\
&=& \sum_{i=1}^{n} A_{vi}^{\rm T}\Lambda^{3*}_i=A^{\rm T}\Lambda^{3*},\label{rcc-opt-tmp1}
\\
0_{m\times p}&=&\sum_{i=1}^{n}[-(Y^*B_{li} - F_{li})B_{li}^{\rm T}+{[I_{m_i}]_{\rm R} \Lambda^{3*}_i} \nonumber\\
&&-\sum_{j=1}^n a_{j,i}(\Lambda^{2*}_i-\Lambda^{2*}_j)]\nonumber\\
&=&
\sum_{i=1}^{n}[-(Y^*B_{li} - F_{li})B_{li}^{\rm T}+{[I_{m_i}]_{\rm R} \Lambda^{3*}_i} ]\nonumber\\
&=&-(Y^*B - F)B^{\rm T}+ \Lambda^{3*},\label{rcc-opt-tmp2}
\end{eqnarray}
where $\Lambda^{3*} = \begin{bmatrix}(\Lambda^{3*}_1)^{\rm T}
\cdots
(\Lambda^{3*}_n)^{\rm T}
\end{bmatrix}^{\rm T}\in\mathbb R^{m\times p}$. It follows from \eqref{rcc-opt-tmp1} and \eqref{rcc-opt-tmp2} that $A^{\rm T}(Y^*B-F)B^{\rm T}=0_{r\times p}$. Recall $AX^*=Y^*$. Eqn. \eqref{opt-cond-c} holds and $X^*$ is a least squares solution to \eqref{AXBeF}.

($ii$) Conversely, suppose that $X^*$ is a least squares solution to \eqref{AXBeF} and $Y^*=AX^*$. We show that $(X_{\rm E}^*,Y_{\rm E}^*)=(1_n\otimes X^*,1_n\otimes Y^*)$  is a solution to problem \eqref{opt-rcc} by proving \eqref{opt-rcc-c}.

Let $\Lambda^{3*} = \begin{bmatrix}(\Lambda^{3*}_1)^{\rm T}
\cdots
(\Lambda^{3*}_n)^{\rm T}
\end{bmatrix}^{\rm T}= (Y^*B-F)B^{\rm T}$. \eqref{opt-cond-c} can be rewritten as
\begin{align*}
&A^{\rm T}\Lambda^{3*}=\sum_{i=1}^n A_{vi}^{\rm T}\Lambda^{3*}_i =0_{r\times p},\\
&\Lambda^{3*}-(Y^*B - F)B^{\rm T}=\sum_{i=1}^{n}[-(Y^*B_{li} - F_{li})B_{li}^{\rm T}+{[I_{m_i}]_{\rm R} \Lambda^{3*}_i} ].
\end{align*}

Because $\ker (L_n)$ and $\mathrm{range} (L_n)$ form an orthogonal decomposition of $\mathbb R^{n}$ by the fundamental theorem of linear algebra \cite{Strang:1993}, where $L_n$ is the Laplacian matrix of $\mathcal G$, there are matrices $\Lambda_i^{1*}\in\mathbb R^{r\times p}$, $\Lambda_i^{2*}\in\mathbb R^{m\times p}$ such that \eqref{opt-rcc-c} holds. It follows from $AX^*=Y^*$ and the KKT optimality condition (Theorem 3.25 of \cite{Ruszczynski:2006}) that
$(X_{\rm E}^*,Y_{\rm E}^*)=(1_n\otimes X^*,1_n\otimes Y^*)$  is a solution to problem \eqref{opt-rcc}.

\section{Proof of Proposition \ref{THM_RCC}}\label{app-THM_RCC}

1) Let $(X_{\rm E}^*,Y_{\rm E}^*,\Lambda^{1*}, \Lambda^{2*}, \Lambda^{3*})$  be any equilibrium of algorithm \eqref{algorithm2} and function $V$ be a positive definite function given by
\begin{eqnarray*}
V(\rm {X}_{\rm E},\rm {Y}_{\rm E},\Lambda^1,\Lambda^2,\Lambda^3) & \triangleq & \frac{1}{2}\sum_{i=1}^n \|X_i-X_i^*\|_{\rm F}^2 \nonumber\\
&&\hspace{-2 cm}+  \frac{1}{2}\sum_{i=1}^n \|Y_i-Y_i^*\|_{\rm F}^2 +  \frac{1}{2}\sum_{i=1}^n \|\Lambda^1_i-\Lambda^{1*}_i\|_{\rm F}^2 \nonumber\\
&&\hspace{-2 cm}+ \frac{1}{2} \sum_{i=1}^n \|\Lambda^2_i-\Lambda^{2*}_i\|_{\rm F}^2 + \frac{1}{2} \sum_{i=1}^n \| \Lambda^3_i-\Lambda^{3*}_i\|_{\rm F}^2.
\end{eqnarray*}

The derivative of function $V$ along the trajectory of algorithm \eqref{algorithm2} is given by
\begin{eqnarray}
\dot V &=&  \sum_{i=1}^n \langle X_i-X_i^*, \dot X_i\rangle_{\rm F} + \sum_{i=1}^n \langle Y_i-Y_i^*, \dot Y_i\rangle_{\rm F} \nonumber\\
&&+ \sum_{i=1}^n \langle \Lambda^1_i-\Lambda^{1*}_i, \dot \Lambda^1_i\rangle_{\rm F} +\sum_{i=1}^n \langle \Lambda^2_i-\Lambda^{2*}_i, \dot \Lambda^2_i\rangle_{\rm F}\nonumber\\
&&+\sum_{i=1}^n \langle \Lambda^3_i-\Lambda^{3*}_i, \dot \Lambda^3_i\rangle_{\rm F} .
\end{eqnarray}

By algorithm \eqref{algorithm2}  and the facts that $A_{vi}X_i^*-Y_i^{vi*}=0_{m_i\times p}$, $X_i^*=X_j^*$, $-A_{vi}^{\rm T}\Lambda^{3*}_i- \sum_{j=1}^n a_{i,j}(\Lambda^{1*}_i-\Lambda^{1*}_j)=0_{r\times q}$,  we have
\begin{align}
\sum_{i=1}^n \langle X_i-X_i^*, \dot X_i\rangle_{\rm F}
&=  -\sum_{i=1}^n  \langle X_i-X_i^*,A_{vi}^{\rm T}(\Lambda^3_i-\Lambda^{3*}_i)\rangle_{\rm F}\nonumber\\
&\hspace{-2.5 cm}\quad-\sum_{i=1}^n  \langle X_i-X_i^*, A_{vi}^{\rm T}(A_{vi}X_i-A_{vi}X_i^*+Y_i^{vi*}-Y_i^{vi})\rangle_{\rm F}\nonumber\\
&\hspace{-2.5 cm}\quad+ \sum_{i=1}^n  \sum_{j=1}^n a_{i,j} \langle X_i,\Lambda^{1*}_i-\Lambda^{1*}_j\rangle_{\rm F}\nonumber\\
&\hspace{-2.5 cm}\quad- \sum_{i=1}^n  \sum_{j=1}^n a_{i,j} \langle X_i,\Lambda^1_i-\Lambda^1_j\rangle_{\rm F}\nonumber\\
&\hspace{-2.5 cm}\quad- \sum_{i=1}^n  \sum_{j=1}^n a_{i,j} \langle X_i,X_i-X_j\rangle_{\rm F}\nonumber\\
&\hspace{-2.5 cm}= -\sum_{i=1}^n\|A_{vi}(X_i-X_i^*)\|_{\rm F}^2  -\frac{1}{2} \sum_{i=1}^n  \sum_{j=1}^n a_{i,j} \|X_i-X_j\|_{\rm F}^2\nonumber\\
&\hspace{-2.5 cm}\quad-\sum_{i=1}^n  \langle X_i-X_i^*,A_{vi}^{\rm T}(\Lambda^3_i-\Lambda^{3*}_i)\rangle_{\rm F}
\nonumber\\
& \hspace{-2.5 cm}\quad+ \sum_{i=1}^n\langle X_i-X_i^*,A_{vi}^{\rm T}(Y_i^{vi}-Y_i^{vi*})\rangle_{\rm F}\nonumber\\
&\hspace{-2.5 cm}\quad
-\sum_{i=1}^n \sum_{j=1}^n a_{i,j} \langle \Lambda^1_i-\Lambda^{1*}_i, X_i-X_j \rangle_{\rm F},
\end{align}
\begin{align}
\sum_{i=1}^n \langle Y_i-Y_i^*, \dot Y_i\rangle_{\rm F}
&= -\sum_{i=1}^n  \langle Y_i-Y_i^*, (Y_i - Y_i^*)B_{li}B_{li}^{\rm T}\rangle_{\rm F}\nonumber\\
&\hspace{-1 cm}\quad + \sum_{i=1}^n  \langle Y_i-Y_i^*,[I_{m_i}]_{\rm R} (\Lambda^3_i-\Lambda^{3*}_i)\rangle_{\rm F}\nonumber\\
&\hspace{-1 cm}\quad+\sum_{i=1}^n  \langle Y_i-Y_i^*, [I_{m_i}]_{\rm R}(A_{vi}X_i-Y_i^{vi})\rangle_{\rm F}
\nonumber\\
&\hspace{-1 cm}\quad
-\sum_{i=1}^n \sum_{j=1}^n a_{i,j} \langle Y_i,Y_i-Y_j\rangle_{\rm F}\nonumber\\
&\hspace{-1 cm}\quad -\sum_{i=1}^n \sum_{j=1}^n a_{i,j} \langle Y_i,\Lambda^2_i-\Lambda^2_j\rangle_{\rm F}\nonumber\\
&\hspace{-1 cm}\quad+\sum_{i=1}^n \sum_{j=1}^n a_{i,j} \langle Y_i,\Lambda^{2*}_i-\Lambda^{2*}_j\rangle_{\rm F}\nonumber\\
&\hspace{-1 cm}= -\sum_{i=1}^n  \| (Y_i - Y_i^*)B_{li}\|_{\rm F}^2 \nonumber\\
&\hspace{-1 cm}\quad+ \sum_{i=1}^n  \langle \Lambda^3_i-\Lambda^{3*}_i, Y_i^{vi}-Y_i^{vi*} \rangle_{\rm F},
\nonumber\\
&\hspace{-1 cm}\quad+\sum_{i=1}^n  \langle Y_i^{vi}-Y_i^{vi*},   A_{vi}(X_i-X_i^*)\rangle_{\rm F}\nonumber\\
&\hspace{-1 cm}\quad- \sum_{i=1}^n  \| Y_i^{vi}-Y_i^{vi*} \|_{\rm F}^2
\nonumber\\
&\hspace{-1 cm}\quad -\frac{1}{2}\sum_{i=1}^n \sum_{j=1}^n a_{i,j} \|Y_i-Y_j\|_{\rm F}^2 \nonumber\\
&\hspace{-1 cm}\quad-\sum_{i=1}^n \sum_{j=1}^n a_{i,j} \langle \Lambda^2_i-\Lambda^{2*}_i, Y_i-Y_j \rangle_{\rm F},
\end{align}
\begin{eqnarray}
\sum_{i=1}^n \langle \Lambda^1_i-\Lambda^{1*}_i, \dot \Lambda^1_i\rangle_{\rm F}   =  \sum_{i=1}^n \sum_{j=1}^n a_{i,j} \langle \Lambda^1_i-\Lambda^{1*}_i, X_i-X_j \rangle_{\rm F},
\end{eqnarray}
\begin{eqnarray}
\sum_{i=1}^n \langle \Lambda^2_i-\Lambda^{2*}_i, \dot \Lambda^2_i\rangle_{\rm F}   =  \sum_{i=1}^n \sum_{j=1}^n a_{i,j} \langle \Lambda^2_i-\Lambda^{2*}_i, Y_i-Y_j \rangle_{\rm F},
\end{eqnarray}
\begin{align}
\sum_{i=1}^n \langle \Lambda^3_i-\Lambda^{3*}_i, \dot \Lambda^3_i\rangle_{\rm F}   = & \sum_{i=1}^n  \langle \Lambda^3_i-\Lambda^{3*}_i, A_{vi}(X_i- X_i^*)\rangle_{\rm F}\nonumber\\
& - \sum_{i=1}^n  \langle \Lambda^3_i-\Lambda^{3*}_i, Y_i^{vi}-Y_i^{vi*} \rangle_{\rm F}.
\end{align}

To sum up,
\begin{align}\label{DV_RCC}
\dot V
&= -\frac{1}{2}\sum_{i=1}^n \sum_{j=1}^n a_{i,j} \|X_i-X_j\|_{\rm F}^2-\sum_{i=1}^n\|A_{vi}X_i-Y_i^{vi}\|_{\rm F}^2 \nonumber\\
& -\sum_{i=1}^n  \| (Y_i - Y_i^*)B_{li}\|_{\rm F}^2  - \frac{1}{2}\sum_{i=1}^n \sum_{j=1}^n a_{i,j} \|Y_i-Y_j\|_{\rm F}^2\leq 0.
\end{align}

Hence, $(X_{\rm E}^*,Y_{\rm E}^*,\Lambda^{1*}, \Lambda^{2*}, \Lambda^{3*})$ is a Lyapunov stable equilibrium of algorithm \eqref{algorithm2}. Because function $V$ is positive definite and radically unbounded. It follows from \eqref{DV_RCC} that a trajectory of algorithm \eqref{algorithm2} is bounded for arbitrary initial condition.

2)  Define the set
\begin{eqnarray*}
\mathcal R &=&\big{\{}(X_{\rm E},Y_{\rm E},\Lambda^1,\Lambda^2,\Lambda^3):\dot V(X_{\rm E},Y_{\rm E},\Lambda^1,\Lambda^2,\Lambda^3)
=0 \big{\}}\\
&\subset &\big{\{}(X_{\rm E},Y_{\rm E},\Lambda^1,\Lambda^2,\Lambda^3):A_{vi}X_i-Y_i^{vi}=0_{m_i\times p},\\
&&\,(Y_i-Y_i^*)B_{li}=0_{m\times q_i},\,X_i=X_j,\,Y_i=Y_j,\\
&&\,i,j\in\{1,\ldots,n\}\big{\}}.
\end{eqnarray*}

Let $\mathcal M$ be the largest invariant subset of $\overline{\mathcal R}$.
It follows from the invariance principle (Theorem 2.41 of \cite{HC:2008}) that
$(X_{\rm E}(t),Y_{\rm E}(t),\Lambda^1(t),\Lambda^2(t),\Lambda^3(t))\rightarrow\mathcal M$ as $t\rightarrow \infty$  and $\mathcal M$ is positive invariant.
Assume that $(\overline X_{\rm E}(t),\overline Y_{\rm E}(t),\overline \Lambda^1(t),\overline \Lambda^2(t),\overline \Lambda^3(t))$ is a trajectory of \eqref{algorithm2} with $$(\overline X_{\rm E}(t),\overline Y_{\rm E}(t),\overline \Lambda^1(t),\overline \Lambda^2(t),\overline \Lambda^3(t))\in \mathcal M$$ for all $t\geq 0$. For all $i\in\{1,\ldots,n\}$, we have
$\dot {\overline\Lambda}_i^1(t)\equiv 0_{r\times q},$ $\dot {\overline \Lambda}_i^2(t)\equiv 0_{m\times q},$ and $\dot {\overline\Lambda}_i^3(t)\equiv 0_{m_i\times p}$ and hence,
\begin{eqnarray*}
\dot {\overline { X}}_i(t)\equiv -A_{vi}^{\rm T}\overline\Lambda^3_i(0)- \sum_{j=1}^n a_{i,j}(\overline\Lambda^1_i(0)-\overline\Lambda^1_j(0)),
\end{eqnarray*}
\begin{eqnarray*}
\dot {\overline { Y}}_i(t) &=& -(\overline Y_i(t)B_{li} - F_{li})B_{li}^{\rm T}+{[I_{m_i}]_{\rm R} \Lambda^3_i}(0)  \\
&&-\sum_{j=1}^n a_{i,j}(\overline\Lambda^2_i(0)-\overline\Lambda^2_j(0))\\
&=& -(\overline Y_i(t) -Y_i^*)B_{li}B_{li}^{\rm T}-(Y_i^*B_{li}- F_{li})B_{li}^{\rm T} \\
 &&+{[I_{m_i}]_{\rm R} \overline\Lambda^3_i}(0)-\sum_{j=1}^n a_{i,j}(\overline\Lambda^2_i(0)-\overline\Lambda^2_j(0))\\
&\equiv & -(Y_i^*B_{li}- F_{li})B_{li}^{\rm T}+{[I_{m_i}]_{\rm R} \overline\Lambda^3_i}(0) \\
 &&-\sum_{j=1}^n a_{i,j}(\overline\Lambda^2_i(0)-\overline\Lambda^2_j(0)).
\end{eqnarray*}

Suppose $\dot {\overline { X}}_i(t)\not=0_{r\times p}$ (or $\dot {\overline { Y}}_i(t)\not=0_{m\times p}$).  Then ${\overline  { X}}_i(t)\rightarrow\infty$ (or ${\overline  { Y}}_i(t)\rightarrow\infty$) as $t\rightarrow\infty$, which contradicts the boundedness of the trajectory. Hence, $\dot {\overline  { X}}_i(t)=0_{r\times p}$, $\dot {\overline { Y}}_i(t)=0_{m\times p}$, and
$\mathcal M\subset \big{\{}(X_{\rm E},Y_{\rm E},\Lambda^1,\Lambda^2,\Lambda^3): \dot {{ X}}_i=0_{r\times p},\,\dot {{ Y}}_i=0_{m\times p},\,\dot\Lambda_i^1=0_{r\times q},\,\dot\Lambda_i^2=0_{m\times q},\, \dot\Lambda_i^3=0_{m_i\times p}\big{\}}$.

Clearly, any point  in $\mathcal M$ is   an equilibrium point of algorithm \eqref{algorithm2}. By  part 1), any point  in $\mathcal M$  is Lyapunov stable. It follows from Lemma \ref{semistability} that    \eqref{algorithm2} is globally convergent  to an  equilibrium.
Due to Proposition \ref{rcc-se} and Lemma \ref{rcc-eq}, $X_i(t)$ converges to a least squares solution to \eqref{AXBeF}. Furthermore, it follows form Lemma \ref{exponential} that the convergence rate of algorithm \eqref{algorithm2} is exponential.

\section{Proof of Proposition \ref{Thm-RRR}}\label{app-Thm-RRR}

1) Let $(X^*,Y_{\rm E}^*,Z^*, \Lambda^{1*},\Lambda^{2*})$   be an equilibrium of algorithm \eqref{algorithm1} and define function $V$ as \eqref{V_RRR}.
The function derivative $\dot V_1(\cdot)$ 
along the trajectory of algorithm \eqref{algorithm1} is
\begin{eqnarray*}
\dot V_1  &=& \sum_{i=1}^n \Big\langle A_{vi}^{\rm T}(A_{vi}Y_i -F_{vi})+\sum_{j=1}^n a_{i,j}(Y_i-Y_j)\\
&&-A_{vi}^{\rm T}(A_{vi}Y_i^* -F_{vi}), \dot Y_i \Big \rangle_{\rm F},\\
&=&\sum_{i=1}^n \Big\langle A_{vi}^{\rm T}(A_{vi}Y_i -F_{vi}) +\frac{1}{n}\Lambda^1_i\\
&&+\sum_{j=1}^n a_{i,j}(Y_i-Y_j) +\sum_{j=1}^n a_{i,j}(\Lambda^2_i-\Lambda^2_j), \dot Y_i \Big\rangle_{\rm F} \\
&&+\sum_{i=1}^n \Big\langle -\frac{1}{n}\Lambda^1_i-\sum_{j=1}^n a_{i,j}(\Lambda^2_i-\Lambda^2_j)\\
&&-A_{vi}^{\rm T}(A_{vi}Y_i^* -F_{vi}), \dot Y_i \Big\rangle_{\rm F}.
\end{eqnarray*}
Note that $-A_{vi}^{\rm T}(A_{vi}Y_i^*-F_{vi}) -\frac{1}{n}\Lambda^{1*}_i -\sum_{j=1}^n a_{i,j}(\Lambda^{2*}_i-\Lambda^{2*}_j)=0_{r\times q}$ because $(X^*,Y_{\rm E}^*,Z^*, \Lambda^{1*},\Lambda^{2*})$   is an equilibrium of algorithm \eqref{algorithm1}. Thus,
\begin{eqnarray*}
\dot V_1  &=& \sum_{i=1}^n \Big\langle A_{vi}^{\rm T}(A_{vi}Y_i -F_{vi}) +\frac{1}{n}\Lambda^1_i+\sum_{j=1}^n a_{i,j}(Y_i-Y_j)\\
&& +\sum_{j=1}^n a_{i,j}(\Lambda^2_i-\Lambda^2_j), \dot Y_i \Big\rangle_{\rm F} -\frac{1}{n}\sum_{i=1}^n \langle\Lambda^1_i-\Lambda^{1*}_i, \dot Y_i \rangle_{\rm F}\\
&&-\sum_{i=1}^n \sum_{j=1}^n a_{i,j} \langle \Lambda^2_i-\Lambda^2_j-\Lambda^{2*}_i+\Lambda^{2*}_j, \dot Y_i \rangle_{\rm F}\\
&=&-\|\dot Y_i\|_{\rm F}^2-\frac{1}{n}\sum_{i=1}^n \langle\Lambda^1_i-\Lambda^{1*}_i, \dot Y_i \rangle_{\rm F}\\
&&-\sum_{i=1}^n \sum_{j=1}^n a_{i,j}\langle \Lambda^2_i-\Lambda^{2*}_i, \dot Y_i- \dot Y_j\rangle_{\rm F}.
\end{eqnarray*}

Following similar steps to prove  part 1) of Proposition \ref{THM_RCC}, we can prove that $\dot V_2$ along the trajectory of algorithm \eqref{algorithm1} is
\begin{eqnarray*}
\dot V_2   &=& -\sum_{i=1}^n  \| A_{vi}(Y_i - Y_i^*)\|_{\rm F}^2-\frac{1}{2}\sum_{i=1}^n \sum_{j=1}^n a_{i,j} \|Y_i-Y_j\|_{\rm F}^2\\
&& +\frac{1}{n}\sum_{i=1}^n \langle\Lambda^1_i-\Lambda^{1*}_i, \dot Y_i \rangle_{\rm F} -\frac{1}{2}\sum_{i=1}^n \sum_{j=1}^n a_{i,j} \|\Lambda^1_i-\Lambda^1_j\|_{\rm F}^2 \nonumber\\
&& - \|\Lambda^1_i B_{vi}^{\rm T}\|^2_{\rm F} +\sum_{i=1}^n \sum_{j=1}^n a_{i,j}\langle \Lambda^2_i-\Lambda^{2*}_i, \dot Y_i -\dot Y_j \rangle_{\rm F}.
\end{eqnarray*}

Hence,
\begin{eqnarray}\label{DV_RRR}
\dot V    &=& -\|\dot Y_i\|_{\rm F}^2-\sum_{i=1}^n  \| A_{vi}(Y_i - Y_i^*)\|_{\rm F}^2- \|\dot X_{li}\|^2_{\rm F}\nonumber\\
&&
-\frac{1}{2}\sum_{i=1}^n \sum_{j=1}^n a_{i,j} \|Y_i-Y_j\|_{\rm F}^2\nonumber\\
&&-\frac{1}{2}\sum_{i=1}^n \sum_{j=1}^n a_{i,j} \|\Lambda^1_i-\Lambda^1_j\|_{\rm F}^2 \leq 0.
\end{eqnarray}

Recall that $V_1(Y_{\rm E})$  is nonnegative for all $Y_{\rm E}\in \mathbb R^{nr\times q}$ due to Lemma \ref{V_rrr_p}.  $V$ is positive definite and radically unbounded, $(X^*,Y_{\rm E}^*,Z^*, \Lambda^{1*},\Lambda^{2*})$  is a Lyapunov stable equilibrium, and furthermore, it follows from \eqref{DV_RRR} that a trajectory of algorithm \eqref{algorithm1} is bounded for arbitrary initial condition.

2) Let
\begin{align*}
\mathcal R =&\big{\{}({ X},Y_{\rm E},Z,\Lambda^1,\Lambda^2):\dot V(X,Y_{\rm E}, Z,\Lambda^1,\Lambda^2)
=0 \big{\}}\\
\subset &\big{\{}({ X},Y_{\rm E},Z,\Lambda^1,\Lambda^2):A_{vi}(Y_i - Y_i^*)=0_{m_i\times q},\,\dot Y_i = 0_{r\times q},\\
&\Lambda^1_i=\Lambda^1_j,\,Y_i=Y_j,\,\dot X_{li}=0_{r\times p_i},\,i,j\in\{1,\ldots,n\}\big{\}}.
\end{align*}

Let $\mathcal M$ be the largest invariant subset of $\overline{\mathcal R}$.
It follows from the invariance principle (Theorem 2.41 of \cite{HC:2008}) that
$({X}(t),Y_{\rm E}(t),Z(t),\Lambda^1(t),\Lambda^2(t))\rightarrow\mathcal M$ as $t\rightarrow \infty$. Note that $\mathcal M$ is invariant. The trajectory $({X}(t),Y_{\rm E}(t),Z(t),\Lambda^1(t),\Lambda^2(t))\in \mathcal M$ for all $t\geq 0$ if $$({X}(0),Y_{\rm E}(0),Z(0),\Lambda^1(0),\Lambda^2(0))\in\mathcal M.$$ Assume $(\overline{X}(t),\overline Y_{\rm E}(t),\overline Z(t),\overline\Lambda^1(t),\overline\Lambda^2(t))\in \mathcal M$ for all $t\geq 0$, $\dot {\overline { X}}_{li}(t) \equiv 0_{r\times p_i}$, $\dot  {\overline { Y}}_{i}(t)\equiv 0_{r\times q}$,
$\dot {\overline Z}_i(t)\equiv 0_{r\times q},$  $\dot {\overline \Lambda}_i^2(t)\equiv 0_{r\times q},$ and hence,
\begin{eqnarray*}
\dot {\overline \Lambda}^1_i(t) = \frac{1}{n}\overline Y_i(0)-\overline X_{li}(0)B_{vi} +\sum_{j=1}^n a_{i,j}(\overline Z_i(0)-\overline Z_j(0)).
\end{eqnarray*}

If $\dot {\overline \Lambda}^1_i(t)\not=0_{r\times q}$, then ${\overline \Lambda}^1_i(t)\rightarrow\infty$  as $t\rightarrow\infty$, which contradicts the boundedness of the trajectory.  Hence, $\dot {\overline \Lambda}^1_i(t)\equiv 0_{r\times q}$ for $i\in\{1,\ldots,n\}$ and $\mathcal M\subset \big{\{}({X},Y_{\rm E},Z,\Lambda^1,\Lambda^2): \dot {{ X}}_{li}=0_{r\times p_i},\,\dot {{ Y}}_i=0_{r\times q},\,\dot {{ Z}}_i=0_{r\times p},\,\dot\Lambda_i^1=0_{r\times q},\,\dot\Lambda_i^2=0_{m\times q}\big{\}}$.

Take any $(\tilde {X},\tilde Y_{\rm E},\tilde Z,\tilde\Lambda^1,\tilde \Lambda^2)\in \mathcal M$. $(\tilde {X},\tilde Y_{\rm E},\tilde Z,\tilde\Lambda^1,\tilde \Lambda^2)$ is clearly an equilibrium point of algorithm \eqref{algorithm1}. It follows from  part 1) that $(\tilde {X},\tilde Y_{\rm E},\tilde Z,\tilde\Lambda^1,\tilde \Lambda^2)$  is Lyapunov stable. Hence, every point in $\mathcal M$ is Lyapunov stable.
By Lemma \ref{semistability}, algorithm (\ref{algorithm1}) is  convergent
to
an equilibrium.
Due to Lemma \ref{exponential}, algorithm (\ref{algorithm1}) converges to an equilibrium exponentially.
According to Proposition \ref{rrr-se} and Lemma \ref{rrr-eq}, $X(t)$ converges to a least squares solution to  \eqref{AXBeF} exponentially.

\section{Proof of Proposition \ref{THM3}}\label{app-THM3}

1) Let functions $V_1$ and $V_2$ be as defined in \eqref{V_CCR}.
Their derivatives along the trajectory of algorithm \eqref{algorithm3} are
\begin{eqnarray*}
\dot V_1  &=& -\sum_{i=1}^n\|\dot U_i\|_{\rm F}^2- \sum_{i=1}^n \langle \Lambda^2_i-\Lambda^{2*}_i, \dot U_i\rangle_{\rm F}\\
&&-\sum_{i=1}^n\|\dot Y_{vi}\|_{\rm F}^2 -\sum_{i=1}^n \langle \Lambda^3_i-\Lambda^{3*}_i, [\dot Y_{vi}]_{\rm R}\rangle_{\rm F},
\end{eqnarray*}
\begin{eqnarray}
\dot V_2    &=&  \sum_{i=1}^n \langle X_{i}-X_{i}^*, \dot X_{i}\rangle_{\rm F} + \sum_{i=1}^n \langle Y_{vi}-Y_{vi}^*, \dot Y_{vi}\rangle_{\rm F} \nonumber\\
&&  + \sum_{i=1}^n \langle Z_i-Z_i^*, \dot Z_i\rangle_{\rm F}+ \sum_{i=1}^n \langle U_i-U_i^*, \dot U_i\rangle_{\rm F}\nonumber\\
&&   + \sum_{i=1}^n \langle W_i-W_i^*, \dot W_i\rangle_{\rm F}+ \sum_{i=1}^n \langle \Lambda^1_i-\Lambda^{1*}_i, \dot \Lambda^1_i\rangle_{\rm F}\nonumber\\
&&  +\sum_{i=1}^n \langle \Lambda^2_i-\Lambda^{2*}_i, \dot \Lambda^2_i\rangle_{\rm F} +\sum_{i=1}^n \langle \Lambda^3_i-\Lambda^{3*}_i, \dot \Lambda^3_i\rangle_{\rm F}.
\end{eqnarray}

It follows from similar proof of Proposition \ref{Thm-RRR} that the derivative of $V$, which is defined in \eqref{V_CCR}, is
\begin{align}
\dot V
 =& - \frac{1}{2}\sum_{i=1}^n  \sum_{j=1}^n a_{i,j}\|X_i-X_j\|^2_{\rm F}
 -\sum_{i=1}^n  \| A_{li}(Y_{vi} - Y_{vi}^*)\|_{\rm F}^2 \nonumber\\
&-\sum_{i=1}^n\|U_i-U_i^*\|^2_{\rm F}   -\frac{1}{2}\sum_{i=1}^n \sum_{j=1}^n a_{i,j}\|\Lambda^2_i-\Lambda^2_j\|_{\rm F}^2-\sum_{i=1}^n\|\dot U_i\|_{\rm F}^2\nonumber\\
& -\frac{1}{2}\sum_{i=1}^n \sum_{j=1}^n a_{i,j}\|\Lambda^3_i-\Lambda^3_j\|_{\rm F}^2 -\sum_{i=1}^n\|\dot Y_{vi}\|_{\rm F}^2\leq 0.\label{DV_CCR}
\end{align}
Recall that $ V_1(X_{\rm E},Y,U)\geq 0 $ by Lemma \ref{V_1_p}.
Function $V$ is clearly positive definite and radically unbounded. Hence, $( X_{\rm E}^*,Y^*,Z^*,U^*,W^*,\Lambda^{1*},\Lambda^{2*},\Lambda^{3*})$ is Lyapunov stable and the trajectory of algorithm \eqref{algorithm3} is bounded for arbitrary initial condition.

2) Take
\begin{eqnarray*}
\mathcal R &=&\big{\{}(X_{\rm E},Y,Z,U,W,\Lambda^{1},\Lambda^{2},\Lambda^{3}):\dot V =0 \big{\}}\\
&\subset &\big{\{}(X_{\rm E},Y,Z,U,W,\Lambda^{1},\Lambda^{2},\Lambda^{3}):X_i=X_j,\,\Lambda^2_i=\Lambda^2_j,\\
&& A_{li}(Y_{vi} - Y_{vi}^*)=0_{m\times q},\,U_i=U_i^*,\,\Lambda^3_i=\Lambda^3_j,\\
&&\dot Y_{vi}=0_{r_i\times q},\,\dot U_i=0_{m\times q},\,i,j\in\{1,\ldots,n\}\big{\}}.
\end{eqnarray*}

Let $\mathcal M$ be the largest invariant subset of $\overline{\mathcal R}$.
It follows from the invariance principle (Theorem 2.41 of \cite{HC:2008}) that
$(X_{\rm E}(t),Y(t),Z(t),U(t),W(t),\Lambda^{1}(t),\Lambda^{2}(t),\Lambda^{3}(t))\rightarrow\mathcal M$ as $t\rightarrow \infty$. Note that $\mathcal M$ is invariant. The trajectory $(X_{\rm E}(t),Y(t),Z(t),U(t),W(t),\Lambda^{1}(t),\Lambda^{2}(t),\Lambda^{3}(t))\in \mathcal M$ for all $t\geq 0$ if $(X_{\rm E}(0),Y(0),Z(0),U(0),W(0),\Lambda^{1}(0),\Lambda^{2}(0),\Lambda^{3}(0))\in\mathcal M$.

Assume $(\overline X_{\rm E}(t),\overline Y(t),\overline Z(t),\overline U(t), \overline W(t),\overline \Lambda^{1}(t),\overline \Lambda^{2}(t),\overline \Lambda^{3}(t))\in \mathcal M\subset\mathcal R$ for all $t\geq 0$. Then  $\dot {\overline \Lambda}_i^1(t)\equiv 0_{r\times q}$,
$\dot {\overline Z}_i(t)\equiv 0_{r\times q},$ $\dot {\overline U}_i(t) =0_{m\times q}$, $\dot {\overline Y}_{vi}=0_{r_i\times q}$, $\dot {\overline W}_i(t)\equiv 0_{m\times q},$ and hence,
\begin{align}
\dot {\overline X}_i(t) =& \overline \Lambda^3_i(t)[B_{li}]_{\rm C}^{\rm T}-\sum_{j=1}^n a_{i,j}(\overline \Lambda^1_i(0)-\overline \Lambda^1_j(0)),\label{a0-1}\\
\dot  {\overline Y}_{vi}(t) =& -A_{li}^{\rm T}(A_{li} \overline Y_{vi}(0) - [F_{vi}]_{\rm R}-\overline U_i(0))\nonumber\\
&-[I_{r_i}]_{\rm C}\overline \Lambda^3_i(t)=0_{r_i\times q},\label{a0-2}\\
\dot  {\overline \Lambda}^2_i(t) =&   \overline U_i(0)- \sum_{j=1}^n a_{i,j}(\overline W_i(0)-\overline W_j(0)),\label{a0-3}\\
\dot{\overline\Lambda}^3_i(t) =& [\overline Y_{vi}]_{\rm R}(0) - \overline X_{i}(t)[B_{li}]_{\rm C}-\sum_{j=1}^n a_{i,j}(\overline Z_i(0)-\overline Z_j(0)),\nonumber\\
& i\in\{1,\ldots,n\}.\label{a0-4}
\end{align}

If $\dot  {\overline \Lambda}^2_i(t)\not=0_{m\times q}$, then ${\overline \Lambda}^2_i(t)\rightarrow\infty$   as $t\rightarrow\infty$, which contradicts the boundedness of the trajectory. Hence, $\dot  {\overline \Lambda}^2_i(t)=0_{m\times q}$.  Moreover,  \eqref{a0-2} and $\overline \Lambda^3_i(t)=\overline \Lambda^3_j(t)$ imply that
$$\overline \Lambda^3_i(t)\equiv\begin{bmatrix}-A_{l1}^{\rm T}(A_{l1} \overline Y_{v1}(0) - [F_{v1}]_{\rm R}-\overline U_1(0))\\
\vdots\\
-A_{ln}^{\rm T}(A_{ln} \overline Y_{vn}(0) - [F_{vn}]_{\rm R}-\overline U_n(0))
\end{bmatrix}$$
 for all $i\in\{1,\ldots,n\}$. Following similar arguments for proving $\dot  {\overline \Lambda}^2_i(t)=0_{m\times q}$, we have $\dot {\overline \Lambda}^3_i(t)\equiv 0_{r\times q}$ and $\dot {\overline X}_i(t)\equiv 0_{r\times p}$ for all $i\in\{1,\ldots,n\}$. To sum up,
\begin{align*}
\mathcal M  \subset  &\big{\{}(X_{\rm E},Y,Z,U,W,\Lambda^{1},\Lambda^{2},\Lambda^{3}):\dot {X}_i\equiv 0_{r\times p},\,\dot Y_{vi} =0_{r_i\times q},\\
&\dot {Z}_i\equiv 0_{r\times q},\,\dot {U}_i \equiv 0_{m\times q},\,\dot {W}_i\equiv 0_{m\times q},\,\dot {\Lambda}_i^1\equiv 0_{r\times q}, \\
& \dot  {\Lambda}^2_i\equiv 0_{m\times q},\,\dot {\Lambda}^3_i\equiv 0_{r\times q},i\in\{1,\ldots,n\}\big{\}},
\end{align*}
and every  point in $\mathcal M$
is an equilibrium point of algorithm \eqref{algorithm3}.

By part 1), every equilibrium point of algorithm \eqref{algorithm3} is Lyapunov stable, and hence, every point in $\mathcal M$ is a Lyapunov stable equilibrium.
By Lemma \ref{semistability}, algorithm (\ref{algorithm3}) is  convergent to an equilibrium.
As a result, the trajectory of  algorithm (\ref{algorithm3}) converges to a Lyapunov stable equilibrium. Furthermore, it follows from Lemma \ref{exponential} that the trajectory of  algorithm (\ref{algorithm3}) converges to an equilibrium exponentially.
By Proposition \ref{ccr-se} and Lemma \ref{ccr-eq}, $X_i(t)$ converges exponentially and $\lim_{t\rightarrow\infty} X_i(t)$ is a least squares solution to equation \eqref{AXBeF} for all $i\in\{1,\ldots,n\}$.

\section{Proof of Proposition \ref{THM-CRR}}\label{app-THM-CRR}

1) Let function $V$  be as defined in \eqref{V_CRR} and $(X^*,Y^*,Z^*,U^*,W^*,\Lambda^{1*},\Lambda^{2*})$  be any equilibrium of \eqref{algorithm4}. Following similar steps to prove part 1) of Proposition \ref{THM3}, we have
\begin{align}\label{DV_4}
\dot V    =& -\sum_{i=1}^n  \| \dot Y_{vi}\|_{\rm F}^2 -\sum_{i=1}^n \|\dot U_{i}\|^2_{\rm F}-\frac{1}{2}\sum_{i=1}^n \sum_{j=1}^n a_{i,j}\|\Lambda^1_i-\Lambda^1_j\|_{\rm F}^2\nonumber\\
&  -\sum_{i=1}^n \| \Lambda^2_i B_{vi}^{\rm T}\|_{\rm F}^2  -\sum_{i=1}^n  \| (A_{li}Y_{vi}-U_i) - (A_{li}Y_{vi}^*-U_i^*)\|_{\rm F}^2\nonumber\\
& -\frac{1}{2}\sum_{i=1}^n \sum_{j=1}^n a_{i,j}\|\Lambda^2_i-\Lambda^2_j\|_{\rm F}^2\leq 0.
\end{align}

By \eqref{DV_4}, $(X^*,Y^*,Z^*,U^*,W^*,\Lambda^{1*},\Lambda^{2*})$  is  Lyapunov  stable and the trajectory of algorithm \eqref{algorithm4} is bounded for arbitrary initial condition.

2) Let
\begin{align*}
\mathcal R =&\big{\{}( X,Y,Z,U,W,\Lambda^{1},\Lambda^{2}):\dot V(X,Y,Z,U,W,\Lambda^{1},\Lambda^{2} )
=0 \big{\}}\\
\subset &\big{\{}(X,Y,Z,U,W,\Lambda^{1},\Lambda^{2}):\Lambda^2_i B_{vi}^{\rm T}=0_{m\times r_i},\,  \dot Y_{vi}=0_{r_i\times q},\\
&\,\dot U_{i}=0_{m\times q},A_{li}Y_{vi}- U_i=A_{li} Y_{vi}^*- U_i^*,\\
&  \Lambda^1_i=\Lambda^1_j,\,\Lambda^2_i=\Lambda^2_j,\,i,j\in\{1,\ldots,n\}\big{\}}.
\end{align*}

Let $\mathcal M$ be the largest invariant subset of $\overline{\mathcal R}$.
It follows from the invariance principle (Theorem 2.41 of \cite{HC:2008}) that
$(X(t),Y(t),Z(t),U(t),W(t),\Lambda^{1}(t),\Lambda^{2}(t))\rightarrow\mathcal M$ as $t\rightarrow \infty$ and $\mathcal M$ is invariant.
Assume $(\overline {X}(t),\overline Y(t),\overline Z(t),\overline U(t), \overline W(t),\overline \Lambda^{1}(t),\overline \Lambda^{2}(t))\in \mathcal M$ for all $t\geq 0$, $\dot {\overline X}_{li}(t)\equiv 0_{m\times q}$, $\dot {\overline Y}_{vi}(t)\equiv 0_{r_i\times q}$,
$\dot {\overline Z}_i(t)\equiv 0_{r\times q},$ $\dot {\overline U}_i(t) \equiv 0_{m\times q}$, $\dot {\overline W}_i(t)\equiv 0_{m\times q},$ and hence,
\begin{eqnarray*}
\dot  {\overline \Lambda}^1_i(t) &=&
\overline U_i(0)- \sum_{j=1}^n a_{i,j}(\overline W_i(0)-\overline W_j(0)),\\
\dot{\overline\Lambda}^2_i(t) &=& [\overline Y_{vi}]_{\rm R}(0) - \overline X_{li}(0)B_{vi}-\sum_{j=1}^n a_{i,j}(\overline Z_i(0)-\overline Z_j(0)).
\end{eqnarray*}
If $\dot  {\overline \Lambda}^1_i(t)\not=0_{m\times q}$ ($\dot  {\overline \Lambda}^2_i(t)\not=0_{r\times q}$), then ${\overline \Lambda}^1_i(t)\rightarrow\infty$  (${\overline \Lambda}^2_i(t)\rightarrow\infty$) as $t\rightarrow\infty$, which contradicts the boundedness of the trajectory. Hence, $\dot  {\overline \Lambda}^1_i(t)=0_{m\times q}$ and $\dot  {\overline \Lambda}^2_i(t)=0_{r\times q}$ for all $i\in\{1,\ldots,n\}$.
To sum up,
\begin{align*}
\mathcal M  \subset  &\big{\{}( X,Y,Z,U,W,\Lambda^{1},\Lambda^{2}):\dot {X}_{li}\equiv 0_{r\times p_i},\,\dot Y_{vi} \equiv 0_{r_i\times q},\\
& \dot {Z}_i\equiv 0_{r\times q},\,\dot {U}_i \equiv 0_{m\times q},\,\dot {W}_i \equiv  0_{m\times q},\,\dot {\Lambda}_i^1\equiv 0_{m\times q}, \\
& \dot  {\Lambda}^2_i\equiv 0_{r\times q},\,i\in\{1,\ldots,n\}\big{\}},
\end{align*}
and every  point in $\mathcal M$
is an equilibrium point of algorithm \eqref{algorithm4}.

By part 1), every equilibrium point of algorithm \eqref{algorithm4} is Lyapunov stable, and hence, every point in $\mathcal M$ is Lyapunov stable. By Lemma \ref{semistability}, algorithm \eqref{algorithm4} is convergent to an equilibrium.
It follows from Proposition \ref{crr-se} and Lemma \ref{crr-eq} that $X(t)$ converges to  a least squares solution to equation \eqref{AXBeF}. In view of Lemma \ref{exponential}, the convergence rate of algorithm \eqref{algorithm4} is exponential.

\end{appendices}

\bibliographystyle{IEEEtran}

\begin{thebibliography}{99}
\bibitem{NOP:2010}
A.~Nedic, A.~Ozdaglar, and P.~A. Parrilo, ``Constrained consensus and
  optimization in multi-agent networks,'' \emph{IEEE Transactions on Automatic
  Control}, vol.~55, no.~4, pp. 922--938, 2010.

\bibitem{YHL:SCL:2015}
P.~Yi, Y.~Hong, and F.~Liu, ``Distributed gradient algorithm for constrained
  optimization with application to load sharing in power systems,''
  \emph{Systems \& Control Letters}, vol.~83, pp. 45--52, 2015.

\bibitem{KCM:A:2015}
S.~S. Kia, J.~Cort\'{e}s, and S.~Mart\'{i}nez, ``Distributed convex
  optimization via continuous-time coordination algorithms with discrete-time
  communication,'' \emph{Automatica}, vol.~55, pp. 254--264, 2015.

\bibitem{QLX:A:2016}
Z.~Qiu, S.~Liu, and L.~Xie, ``Distributed constrained optimal consensus of
  multi-agent systems,'' \emph{Automatica}, vol.~68, pp. 209--215, 2016.

\bibitem{YHX:Cybernetics:2016}
D.~Yuan, D.~W.~C. Ho, and S.~Xu, ``Regularized primal-dual subgradient method
  for distributed constrained optimization,'' \emph{IEEE Transaction on
  Cybernetics}, vol.~46, no.~9, pp. 2109--2118, 2016.

\bibitem{SJ:2012}
G.~Shi and K.~H. Johansson, ``Randomized optimal consensus of multi-agent
  systems,'' \emph{Automatica}, vol.~48, no.~12, pp. 3018--3030, 2012.

\bibitem{ZYH:2017}
X.~Zeng, P.~Yi, and Y.~Hong, ``Distributed continuous-time algorithm for
  constrained convex optimizations via nonsmooth analysis approach,''
  \emph{IEEE Transactions on Automatic Control}, vol.~62, no.~10, pp.
  5227--5233, 2017.

\bibitem{WE:2011}
J.~Wang and N.~Elia, ``Control approach to distributed optimization,'' in
  \emph{Proc. Allerton Conf. Commun., Control and Comput.}, Monticello, IL,
  2010, pp. 557--561.

\bibitem{GC:2014}
B.~Gharesifard and J.~Cort\'{e}s, ``Distributed continuous-time convex
  optimization on weight-balanced digraphs,'' \emph{IEEE Transactions on
  Automatic Control}, vol.~59, pp. 781--786, 2014.

\bibitem{LW:TAC:2015}
Q.~Liu and J.~Wang, ``A second-order multi-agent network for bound-constrained
  distributed optimization,'' \emph{IEEE Transaction on Automatic Control},
  vol.~60, no.~12, pp. 3310--3315, 2015.

\bibitem{SA:TAC:2015}
G.~Shi, B.~D.~O. Anderson, and U.~Helmke, ``Network flows that solve linear
  equations,'' \emph{IEEE Transactions on Automatic Control}, vol.~62, no.~6,
  pp. 2659--2674, 2017.

\bibitem{MLM:TAC:2015}
S.~Mou, J.~Liu, and A.~S. Morse, ``A distributed algorithm for solving a linear
  algebraic equation,'' \emph{IEEE Transactions on Automatic Control}, vol.~60,
  no.~11, pp. 2863--2878, 2015.

\bibitem{WFM:ACC:2016}
L.~Wang, D.~Fullmer, and A.~S. Morse, ``A distributed algorithm with an
  arbitrary initialization for solving a linear algebraic equation,'' in
  \emph{American Control Conference}, Boston, MA, USA, 2016, pp. 1078 -- 1081.

\bibitem{SA:ACC:2016}
G.~Shi and B.~D.~O. Anderson, ``Distributed network flows solving linear
  algebraic equations,'' in \emph{American Control Conference}, Boston, MA,
  USA, 2016, pp. 2864 -- 2869.

\bibitem{AMMH:Arxiv:2015}
B.~Anderson, S.~Mou, A.~Morse, and U.~Helmke, ``Decentralized gradient
  algorithm for solution of a linear equation,'' \emph{Numerical Algebra,
  Control and Optimization}, vol.~6, no.~3, pp. 319--328, 2016.

\bibitem{LLAS:Arxiv:2017}
Y.~Liu, C.~Lageman, B.~Anderson, and G.~Shi, ``An {A}rrow-{H}urwicz-{U}zawa
  type flow as least squares solver for network linear equations,''
  arXiv:1701.03908v1.

\bibitem{LMM:TAC:2017}
J.~Liu, S.~Mou, and A.~S. Morse, ``Asynchronous distributed algorithms for
  solving linear algebraic equations,'' \emph{IEEE Transactions on Automatic
  Control}, vol.~63, no.~2, pp. 372--385, 2017.

\bibitem{LMNM:Automatica:2017}
J.~Liu, A.~S. Morse, A.~Nedic, and T.~Basar, ``Exponential convergence of a
  distributed algorithm for solving linear algebraic equations,''
  \emph{Automatica}, vol.~83, pp. 37--46, 2017.

\bibitem{CZH:CCC:2017}
K.~Cao, X.~Zeng, and Y.~Hong, ``Continuous-time distributed algorithms for
  solving linear algebraic equation,'' in \emph{The 36th Chinese Control
  Conference}, Dalian, China, 2017.

\bibitem{BIG:1980}
A.~Ben-Isreal and T.~N.~E. Greville, \emph{Generalized Inverses: Theory and
  Applications}.\hskip 1em plus 0.5em minus 0.4em\relax New York:
  Springer-Verlag, 2003.

\bibitem{RM:1971}
C.~R. Rao and S.~K. Mitra, \emph{Generalized Inverse of Matrices and Its
  Applications}.\hskip 1em plus 0.5em minus 0.4em\relax John Wiley \& Sons Inc,
  1972.

\bibitem{BK:1980}
J.~K. Baksalary and R.~Kala, ``The matrix equation {$AXB+CYD=E$},''
  \emph{Linear Algebra and its Applications}, vol.~30, pp. 141--147, 1980.

\bibitem{TW:2013}
Y.~Tian and H.~Wang, ``Relations between least-squares and least-rank solutions
  of the matrix equation {$AXB=C$},'' \emph{Applied Mathematics and
  Computation}, vol. 219, no.~20, pp. 10\,293--10\,301, 2013.

\bibitem{Woodgate:1996}
K.~G. Woodgate, ``Least-squares solution of {$F=PG$} over positive semidefinite
  symmetric {$P$},'' \emph{Linear Algebra and its Applications}, vol. 145, pp.
  171--190, 1996.

\bibitem{WC:1996}
L.~Wu and B.~Cain, ``The re-nonnegative definite solutions to the matrix
  inverse problem {$AX = B$},'' \emph{Linear Algebra and its Applications},
  vol. 236, pp. 137--146, 1996.

\bibitem{MHZ:2005}
C.~J. Meng, X.~Y. Hu, and L.~Zhang, ``The skew-symmetric orthogonal solutions
  of the matrix equation {$AX=B$},'' \emph{Linear Algebra and its
  Applications}, vol. 402, pp. 303--318, 2005.

\bibitem{HHB:TAC:2009}
Q.~Hui, W.~M. Haddad, and S.~P. Bhat, ``Semistability, finite-time stability,
  differential inclusions, and discontinuous dynamical systems having a
  continuum of equilibria,'' \emph{IEEE Transactions on Automatic Control},
  vol.~54, no.~10, pp. 2465--2470, 2009.

\bibitem{Antipin:2003}
A.~S. Antipin, ``Feedback-controlled saddle gradient processes,''
  \emph{Automation and Remote Control}, vol.~55, no.~3, pp. 311--320, 2003.

\bibitem{GR2001}
C.~Godsil and G.~F. Royle, \emph{Algebraic Graph Theory}.\hskip 1em plus 0.5em
  minus 0.4em\relax New York: Springer-Verlag, 2001.

\bibitem{Bernstein:2009}
D.~S. Bernstein, \emph{Matrix Mathematics: Theory, Facts, and Formulas, 2nd
  edition}.\hskip 1em plus 0.5em minus 0.4em\relax Princeton, NJ: Princeton
  University Press, 2009.

\bibitem{DLD:2008}
F.~Ding, P.~X. Liu, and J.~Ding, ``Iterative solutions of the generalized
  {S}ylvester matrix equations by using the hierarchical identification
  principle,'' \emph{Applied Mathematics and Computation}, vol. 197, pp.
  41--50, 2008.

\bibitem{YonggeTian:2009}
Y.~Tian, ``Some properties of submatrices in a solution to the matrix equation
  {$AXB=C$} with applications,'' \emph{IEEE Transactions on Control of Network
  Systems}, vol. 346, pp. 557--569, 2009.

\bibitem{PHZ:2005}
Y.~Peng, X.~Hu, and L.~Zhang, ``An iteration method for the symmetric solutions
  and the optimal approximation solution of the matrix equation {$AXB=C$},''
  \emph{Applied Mathematics and Computation}, vol. 160, pp. 763--777, 2005.

\bibitem{HPZ:2006}
J.~Hou, Z.~Peng, and X.~Zhang, ``An iterative method for the least squares
  symmetric solution of matrix equation {$AXB = C$},'' \emph{Numerical
  Algorithms}, vol.~42, no.~2, pp. 181--192, 2006.

\bibitem{ZLGZ:2011}
F.~Zhang, Y.~Li, W.~Guo, and J.~Zhao, ``Least squares solutions with special
  structure to the linear matrix equation {$AXB = C$},'' \emph{Applied
  Mathematics and Computation}, vol. 217, pp. 10\,049--10\,057, 2011.

\bibitem{Ruszczynski:2006}
A.~Ruszczynski, \emph{Nonlinear Optimization}.\hskip 1em plus 0.5em minus
  0.4em\relax Princeton, New Jersey: Princeton University Press, 2006.

\bibitem{WZMC:arx:2017}
X.~Wang, J.~Zhou, S.~Mou, and M.~J. Corless, ``A distributed algorithm for
  least square solutions of linear equations,'' arXiv:1709.10157v1.

\bibitem{Strang:1993}
G.~Strang, ``The fundamental theorem of linear algebra,'' \emph{American
  Mathematical Monthly}, vol. 100, no.~9, pp. 848--855, 1993.

\bibitem{YHL:Automatica:2016}
P.~Yi, Y.~Hong, and F.~Liu, ``Initialization-free distributed algorithms for
  optimal resource allocation with feasibility constraints and application to
  economic dispatch of power systems,'' \emph{Automatica}, vol.~74, pp.
  259--269, 2016.

\bibitem{HC:2008}
W.~M. Haddad and V.~Chellaboina, \emph{Nonlinear Dynamical Systems and Control:
  A {L}yapunov-Based Approach}.\hskip 1em plus 0.5em minus 0.4em\relax
  Princeton, NJ: Princeton Univ. Press, 2008.
\end{thebibliography}

\end{document}